\documentclass{lms}

\usepackage{amsmath,amsfonts,amscd,pifont} 
\usepackage[dvips]{graphicx}
\usepackage{pb-diagram} 

\newtheorem{theorem}{Theorem}[section] 
\newtheorem{lemma}[theorem]{Lemma}     
\newtheorem{corollary}[theorem]{Corollary}

\newnumbered{assertion}{Assertion}    
\newnumbered{conjecture}{Conjecture}  
\newnumbered{definition}{Definition}
\newnumbered{hypothesis}{Hypothesis}
\newnumbered{remark}{Remark}
\newnumbered{note}{Note}
\newnumbered{observation}{Observation}
\newnumbered{problem}{Problem}
\newnumbered{question}{Question}
\newnumbered{algorithm}{Algorithm}
\newnumbered{example}{Example}
\newunnumbered{notation}{Notation} 

\newcommand{\mding}[1]{\mbox{\ding{#1}}}

\newcommand{\BB}[1]{\ensuremath{\mathbb{#1}}}
\newcommand{\Q}{\ensuremath{\BB{Q}}}
\newcommand{\R}{\ensuremath{\BB{R}}}
\newcommand{\Z}{\ensuremath{\BB{Z}}}
\newcommand{\C}{\ensuremath{\BB{C}}}

\newcommand{\bs}{\ensuremath{\boldsymbol}}
\newcommand{\mf}{\ensuremath{\mathfrak}}

\newcommand{\la}{\ensuremath{\langle}}
\newcommand{\ra}{\ensuremath{\rangle}}

\newcommand{\nrm}[1]{\ensuremath{\mathnormal{\lowercase{#1}}}}

\DeclareMathOperator{\sgn}{sgn}
\DeclareMathOperator{\vol}{vol}

\DeclareMathOperator{\Jac}{Jac}
\DeclareMathOperator{\Pf}{Pf}

\title[The Range of Multiplicative Functions on \mbox{\C[\nrm{x}], \R[\nrm{x}]}
and \mbox{\Z[\nrm{x}]}]
 {The Range of Multiplicative Functions on $\C[x], \R[x]$ and $\Z[x]$} 

\author{Christopher D. Sinclair}

\classno{11B83 (primary), 11J71, 37A45, 60G10 (secondary)}

\begin{document}
\maketitle

\begin{abstract}
  Mahler's measure is generalized to create the class of {\it
    multiplicative distance functions}.  These functions measure the
  complexity of polynomials based on the location of their zeros in
  the complex plane.  Following work of \mbox{S.-J.~Chern} and
  J.~Vaaler in \cite{chern-vaaler}, we associate to each
  multiplicative distance function two families of analytic functions
  which encode information about its range on $\C[x]$ and $\R[x]$.
  These {\it moment functions} are Mellin transforms of distribution
  functions associated to the multiplicative distance function and
  demonstrate a great deal of arithmetic structure.  For instance, we
  show that the moment function associated to Mahler's measure
  restricted to real reciprocal polynomials of degree $2N$ has an
  analytic continuation to rational functions with rational
  coefficients, simple poles at integers between $-N$ and $N$, and a
  zero of multiplicity $2N$ at the origin.  This discovery leads to
  asymptotic estimates for the number of reciprocal integer
  polynomials of fixed degree with Mahler measure less than $T$ as $T
  \rightarrow \infty$.  To explain the structure of this moment
  functions we show that the real moment functions of a multiplicative
  distance function can be written as Pfaffians of antisymmetric
  matrices formed from a skew-symmetric bilinear form associated to
  the multiplicative distance function.
\end{abstract}

\section{Introduction}
This manuscript is concerned with measures of complexity of
polynomials which respect both the algebraic structure and topology of
$\C[x]$ (as generated by all open sets in all finite dimensional
subspaces of $\C[x]$).  As such, we are interested in functions from
$\C[x]$ to the non-negative reals which are continuous (as a function
on coefficient vectors) on all finite dimensional subspaces of $\C[x]$
and behave nicely with respect to multiplication and scalar
multiplication.  The most important requirement of these functions is
that they be multiplicative.  As we shall see, multiplicativity is a
very strong condition which allows for many interesting theorems.

The following axiomatization suggests itself:  A function $\Phi: \C[x]
\rightarrow [0, \infty)$ 
will be called a {\it multiplicative distance function} if   
\begin{enumerate}
\item[A1.]\label{item:1} $\Phi$ is continuous,
\end{enumerate}
and for all $f,g \in \C[x]$ and $w \in \C$,
\begin{enumerate}
\item[A2.]\label{item:2} $\Phi$ is positive definite:  $\Phi(f) = 0$ if and only if
  $f$ is identically zero,
\item[A3.]\label{item:3} $\Phi$ is absolutely homogeneous:  $\Phi(w f) = |w|\Phi(f)$, and
\item[A4.]\label{item:4} $\Phi$ is multiplicative: $\Phi(f g) =
  \Phi(f) \Phi(g)$. 
\end{enumerate}
The nomenclature stems from the fact that multiplicative distance
functions restricted to finite dimensional subspaces of $\C[x]$ are
distance functions in the sense of the geometry of numbers.  We will
refer to $\Phi(f)$ as the $\Phi$-distance of $f$ to the origin or
simply the {\it distance} of $f$.

It is easily seen that $\Phi$ is uniquely determined by its action on
monic linear polynomials.  That is, if 
\begin{equation}
\label{eq:2}
f(x) = a \prod_{n=1}^N (x - \gamma_n),
\end{equation}
then there exists a function $\phi: \C \rightarrow (0, \infty)$ such
that
\[
\Phi(f) = |a| \prod_{n=1}^N \phi(\gamma_n).
\]
The function $\phi$ will be known as the {\it root function} of $\Phi$.
The best known example of a multiplicative distance function is
Mahler's measure, denoted by $\mu$, and defined by
\[
\mu(f) = |a| \prod_{n=1}^N \max\{1, |\gamma_n|\}.
\]
From the definition of $\mu$ it is clear that Mahler's measure
satisfies Axioms A2, A3 and A4.  It is less clear that Mahler's
measure is continuous on $\C[x]$, but this was proved by K.~Mahler in
1961 \cite{mahler}.

Mahler's measure can be extended to the algebra of Laurent polynomials
by using $1 = \mu(1) = \mu(x x^{-1})$ to write $\mu(x^{-1}) =
\mu(x)^{-1} = 1.$ If $\Phi(x^n) = 1$ for every $n \geq 0$ (or what
amounts to the same thing, $\phi(0) = 1$) then we shall say $\Phi$ is
{\it shift invariant}.  Shift invariant multiplicative distance
functions can be naturally extended to the algebra of Laurent
polynomials by setting $\phi(x^{-1}) = 1$.

The continuity of $\Phi$ controls the asymptotic behavior of $\phi$.
And in fact, this asymptotic condition produces a classification of
multiplicative distance functions.
\begin{theorem}
\label{thm:1}
Suppose that $\Phi$ is a multiplicative distance function with root
function $\phi$.  Then,
\[
\phi(\gamma) \sim |\gamma| \qquad \mbox{as} \qquad |\gamma| \rightarrow
\infty. 
\]
Conversely, if $\psi: \C \rightarrow (0, \infty)$ is a continuous
function such that $\psi(\gamma) \sim |\gamma|$ as $|\gamma|
\rightarrow \infty$, then $\psi$ is the root function of a
multiplicative distance function.
\end{theorem}

\subsection{Examples of Multiplicative Distance Functions}

Theorem \ref{thm:1} gives us a way of producing examples of
multiplicative distance functions, and in this section we will
introduce another method for constructing multiplicative distance
functions.

A Laurent polynomial $g(x) \in \C[x, 1/x]$ is said to be {\it reciprocal}
if $g(1/x) = g(x)$, and the algebra of reciprocal Laurent polynomials
is given by $\C[x + 1/x]$.  Clearly the algebra of reciprocal Laurent
polynomials is a subalgebra of $\C[x, 1/x]$ and hence we may speak of
the Mahler measure of a reciprocal Laurent polynomial.  We define {\it
  the reciprocal Mahler's measure}, $\rho$, of $f \in \C[x]$ to be the
Mahler's measure of the reciprocal Laurent polynomial $f(x + 1/x)$.
That is, $\rho(f) = \mu(f(x + 1/x))$.  It follows that the root
function of $\rho$ is given by
\[
\gamma \mapsto \mu(x + 1/x - \gamma) = \max\left\{1, \left|\frac{\gamma + \sqrt{\gamma^2 -
        4}}{2}\right|\right\} \max\left\{1, \left|\frac{\gamma -
      \sqrt{\gamma^2 - 4}}{2}\right| \right\}. 
\]
This definition is independent of the branch of the square root used
and is easily seen to satisfy the conditions of Theorem~\ref{thm:1}.

The procedure used to create the reciprocal Mahler's measure may be
repeated to create multiplicative distance functions formed from
Mahler's measure restricted to other subalgebras of $\C[x, 1/x]$.  In
particular, if $G(x) \in \C[x, 1/x]$ is a fixed Laurent polynomial we
may create a multiplicative distance function by considering Mahler's
measure restricted to $\C[G(x)] \subset \C[x,1/x]$.  Thus we define
$G^{\ast}\mu: \C[x] \rightarrow [0, \infty)$ by $G^{\ast}\mu(f) =
\mu(f \circ G)$.  It can be verified that $G^{\ast}\mu$ satisfies all
the axioms of a multiplicative distance function.  The notation for
$G^{\ast}\mu$ stems from the fact that if we view $G$ as the natural
map $\C[x] \rightarrow \C[G(x)] $ then $G^{\ast} \mu$ is the {\it
  pullback} of $\mu$ through $G$.  That is, $G^{\ast} \mu$ is the map
which makes the following diagram commute.
\[
\begin{diagram}
\node{\C[x]} \arrow{e,t}{G} \arrow{s,l}{G^{\ast} \mu} 
\node{\C[G(x)]} \arrow{sw,b}{\mu} \\
\node{[0,\infty)}
\end{diagram}
\]

Given $0 \leq t \leq 1$ we define the $t$-reciprocal Mahler's measure, $\mu_t$,
to be the pullback of $\mu$ through the Laurent polynomial $x + t/x$.
In this context $\mu_1$ is the reciprocal Mahler's measure and $\mu_0$ is
Mahler's measure.  Thus, as $t$ varies from $0$ to $1$ we have
a `path' of multiplicative distance functions whose end points are $\mu$
and $\rho$.

\subsection{Potentials and Jensen's Formula}
The prototypical multiplicative distance function, Mahler's measure,
satisfies an important integral identity.  If $f$ is given as in (\ref{eq:2})
then Jensen's formula implies that
\begin{equation}
\label{eq:3}
\mu(f) = |a| \prod_{n=1}^N \max\{1, |\gamma_n|\} =
\exp\left\{\frac{1}{2\pi} \int_0^{2 \pi} \log\left|f(e^{i
      \theta})\right| \, d\theta 
\right\}.
\end{equation}
The right hand side of this equation is an example of an equilibrium
potential.  By generalizing the right hand side of (\ref{eq:3}) we
may produce examples of multiplicative distance functions which are
associated to compact subsets of $\C$.  Multiplicative distance
functions of this sort were considered from the standpoint of
equidistribution by R.~Rumely in \cite{rumely}.

Let $K$ be a compact subset of $\C$ and let $\nu$ be a probability
measure whose support is contained in $K$.  The {\it potential} of
$\nu$ is defined to be the function $p_{\nu}: \C \rightarrow
[0,\infty)$ specified by
$$
p_{\nu}(\gamma) = \exp\left\{\int_{K} \log|z - \gamma| \, d\nu(z)  \right\}.
$$It is a fundamental result of potential theory that $p_{\nu}$ is upper
semicontinuous.  Moreover $p_{\nu}(\gamma) \sim |\gamma|$ as $|\gamma|
\rightarrow \infty$ and thus if
$p_{\nu}$ is in fact continuous then it is the root function of a
multiplicative distance function.

If we denote the set of probability measures whose support lies in
$K$ by $M(K)$, then under fairly mild conditions on $K$ there is a unique
probability measure $\nu_K \in M(K)$ which minimizes
\[
I(\nu) = -\int_{K} \log| p_{\nu}(\gamma)| \, d\nu(\gamma) \qquad
\mbox{over all} \qquad \nu
\in M(K).
\]
For instance, the minimizing measure is unique if there exists
at least one $\nu \in M(K)$ with $I(\nu) < \infty$.
When $\nu_K$ exists this measure is known as the {\it equilibrium
  measure} of $K$ and the quantity $c(K) = e^{-I(\nu_K)}$ is known as
the {\it capacity} of $K$.  We will denote the potential of $\nu_K$
simply by $p_K$.  This potential is called the {\it equilibrium
  potential} of $K$.  

If $K$ is regular with respect to the Dirichlet problem then $K$ has
positive capacity and $p_K$ is continuous.  In this situation $p_K$ is
the root function of a multiplicative distance function which will be
denoted $P_K$. For instance, if $K$ is a simply connected compact
subset of $\C$ which does not consist of a single point then $K$ is
regular with respect to the Dirichlet problem and we may speak about
the multiplicative distance function $P_K$.  Explicitly,
\[
P_K(f) = |a| \prod_{n=1}^N p_K(\gamma_n) = \exp\left\{\int_K
  \log|f(z)| \, d\nu_K(z) \right\}.
\]
For example, Mahler's measure can be represented as $P_{D}$ where $D$ is the
closed unit disk.  

It is a well known fact of potential theory that $p_K(\gamma) \geq
c(K)$ with equality if and only if $\gamma \in K$.  Of particular
importance are multiplicative distance functions associated to simply
connected compact sets of capacity 1.  In this situation if $K$
contains the origin then $p_K(0)=1$ and hence $P_K$ is shift invariant.

As the next theorem demonstrates, there is a strong connection between
multiplicative distance functions formed from certain compact sets $K$
and those formed by the pullback of Mahler's measure by certain rational
functions.  

\begin{theorem}
\label{thm:6}
Let $q(x) \in \C[x]$ be a monic polynomials of degree $M$ and define
$G(x) = q(x)/x^{M-1}.$ If $G(x)$ is a conformal map from $\C \setminus
D$ onto its image, then $G^{\ast} \mu = P_K$ where $K$ is the
complement in $\C$ of $G(\C \setminus D)$.
\end{theorem}
Theorem~\ref{thm:6} is well-known to experts in potential theory
(though perhaps not in the language used in this manuscript).
From the definitions of $G^{\ast} \mu$ and $P_K$, the equation
$G^{\ast} \mu = P_K$ may be thought of as an analog of Jensen's
formula.

\begin{corollary}
Let $0 \leq t < 1$ and define $E_t \subset \C$ to be the compact set
given by
\[
E_t = \left\{x + i y: \frac{x^2}{(1 + t)^2} + \frac{y^2}{(1 - t)^2}
  \leq 1 \right\},
\]
and define 
\[E_1 = \{ x + i y :
y = 0, x \in [-2,2] \}.
\]
Then, for any $t \in [0,1]$, the $t$-reciprocal Mahler's measure,
$\mu_t$, is equal to $P_{E_t}$. 
\end{corollary}

As $t$ varies from $0$ to $1$, $E_t$ deforms from the unit disk
through a series of regions bounded by ellipses to the degenerate
ellipse given by the interval $[-2,2]$ on the real axis.  All of these
compact sets have capacity 1.  Thus our `path' of multiplicative
distance functions formed from the pullback of $\mu$ through $x + t/x$
as $t$ varies from $0$ to $1$ can also be thought of as a `path' of
shift invariant multiplicative distance functions formed from the
family of ellipses $E_t$ as $t$ ranges over the same values.

\subsection{Star Bodies and Distribution Functions}
By identifying each polynomial of degree $N$ with its vector of
coefficients, the set of polynomials in $\C[x]$ of
degree $N$ may be identified with the vector space $\C^{N+1}$.   To each
$\mathbf{a} \in \C^{N+1}$ we define the polynomial $\mathbf{a}(x)$ by
\[
\mathbf{a}(x) = \sum_{n=1}^{N+1} a_n x^{N+1-n}.
\]
We may regard $\Phi$ as a function on $\C^{N+1}$ by setting
$\Phi(\mathbf{a}) = \Phi(\mathbf{a}(x))$.  As such $\Phi$ satisfies
all the axioms of a vector norm except the triangle inequality.  The
`unit ball' of $\Phi$ is thus not convex.  That is, the set
\[
\mathcal{V}_N(\Phi) = \{ \mathbf{a} \in \C^{N+1} : \Phi(\mathbf{a})
\leq 1 \}
\]
is a symmetric star body about the origin which will be referred to as
the degree $N$ complex {\it unit star body} of $\Phi$.  Similarly the
degree $N$ {\it real} unit star body is defined to be the set
\[
\mathcal{U}_N(\Phi) = \{ \mathbf{a} \in \R^{N+1} : \Phi(\mathbf{a})
\leq 1 \}.
\]
The absolute homogeneity of $\Phi$ implies that the set of polynomials
of degree $N$ in $\C[x]$ with distance bounded by $T > 0$ is 
the dilated star body $T \mathcal{V}_N$.

As a first application of the theory of multiplicative distance
functions, S-J.~Chern and J.~Vaaler devised a procedure for
determining the volume (Lebesgue measure) of $\mathcal{U}_N(\mu)$ and
then used this to give the main term in an asymptotic estimate for the
number of polynomials in $\Z[x]$ with degree at most $N$ and Mahler
measure bounded by $T$ as $T \rightarrow \infty$ \cite{chern-vaaler}.
Their idea is more generally valid, and we will give similar estimates
for the reciprocal Mahler's measure.
\begin{theorem}
\label{thm:2}
Let $\Phi$ be a multiplicative distance functions.  Then, as $T
\rightarrow \infty$,
\[
\#\left\{\mathbf{a} \in \Z^{N+1} : \Phi(\mathbf{a}) \leq T \right\} =
\vol(\mathcal{U}_N(\Phi)) T^{N+1} + O(T^N). 
\]
\end{theorem}
\begin{proof}
See \cite[Ch. VI, $\S$2]{lang} or \cite[$\S$12]{chern-vaaler}
\end{proof}

In order to determine the volumes of $\mathcal{U}_N(\mu)$ and
$\mathcal{V}_N(\mu)$, Chern and Vaaler introduced two families of
analytic functions which encode information about the range of values
of $\mu$ restricted to polynomials with real and complex coefficients.
Their techniques generalize to other multiplicative distance functions
and the analogous analytic functions demonstrate a great deal of
structure which can be used to learn information about the range of
values of a multiplicative distance function.

We define the degree $N$ {\it monic restriction} of $\Phi$ 
to be the function $\widetilde{\Phi}: \C^N \rightarrow (0,\infty)$ given by 
\[
\widetilde{\Phi}(\mathbf{b}) = \Phi\left(x^N + \sum_{n=1}^N b_n
x^{N-n} \right). 
\]
That is, $\widetilde{\Phi}$ is simply $\Phi$ restricted to the set of
(non-leading) coefficient vectors of monic polynomials of degree $N$. 
We use $\lambda_{N}$ and $\lambda_{2N}$ to denote Lebesgue measure on Borel
subsets of $\R^N$ and $\C^N$ (respectively) and define the
distribution functions $f_N, h_N: [0, \infty) \rightarrow [0,\infty)$
by
\[
f_N(\Phi; \xi) = \lambda_{N}\left\{\mathbf{b} \in \R^N:
\widetilde{\Phi}(\mathbf{b}) \leq \xi \right\},
\]
and
\[
h_N(\Phi; \xi) = \lambda_{2N}\left\{\mathbf{b} \in \C^N:
\widetilde{\Phi}(\mathbf{b}) \leq \xi \right\}.
\]
By identifying $\R^N$ with the set of monic coefficient
vectors in $\R^{N+1}$, $f_N(\Phi; \xi)$ is simply the volume of the intersection of the dilated star body $\xi 
\mathcal{U}_N(\Phi)$ with $\R^N$.  In this way $f_N$
encodes information about the range of values $\Phi$ takes on monic
polynomials of degree $N$ in $\R[x]$.  Similarly $h_N(\Phi; \xi)$
encodes information about the range of values $\Phi$ takes on monic
polynomials of degree $N$ in $\C[x]$.  For instance, the volumes of
$\mathcal{U}_N(\Phi)$ and $\mathcal{V}_N(\Phi)$ can be discovered from
$f_N(\Phi; \xi)$ and $h_N(\Phi; \xi)$.
\begin{theorem}
\label{thm:9}
The supports of $f_N(\Phi; \xi)$ and $h_N(\Phi; \xi)$ are bounded away
from $0$, and as $\xi
\rightarrow \infty$,
\[
f_N(\Phi; \xi) = O(\xi^N) \qquad \mbox{and} \qquad h_N(\Phi; \xi) = O(\xi^{2N}).\]
Moreover,
\[
\lim_{\xi \rightarrow \infty} \frac{f_N(\Phi; \xi)}{\xi^{N}} =
\lambda_N(\mathcal{U}_{N-1}(\Phi)) \qquad \mbox{and} \qquad 
\lim_{\xi \rightarrow \infty} \frac{h_N(\Phi; \xi)}{\xi^{2N}} =
\lambda_{2N}(\mathcal{V}_{N-1}(\Phi)).
\]
\end{theorem}

The Mellin transform of these functions is then given by
\[
\widehat{f_N}(\Phi; s) = \int_0^{\infty} \xi^{-s} f_N(\xi) \,
\frac{d\xi}{\xi} \qquad \mbox{and} \qquad 
\widehat{h_N}(\Phi; s) = \int_0^{\infty} \xi^{-s} h_N(\xi) \,
\frac{d\xi}{\xi},
\]
where $s$ is a complex variable.  From the asymptotic formulae for
$f_N(\Phi; \xi)$ and $h_N(\Phi; \xi)$ it is easy to establish that the
integral defining $\widehat{f_N}(\Phi; s)$ converges when $\Re(s) >
N$, and the integral defining $\widehat{h_N}(\Phi; s)$ converges when
$\Re(s) > 2N$.  Moreover, by Morera's Theorem $\widehat{f_N}$ and
$\widehat{h_N}$ are analytic functions in their respective domains of
convergence.  These analytic functions encode information
about the range of values $\Phi$ takes on monic polynomials of degree
$N$ in $\R[x]$ and $\C[x]$ respectively.  For instance, the volume of
$\mathcal{U}_N(\Phi)$ also appears as a special value of
$\widehat{f_N}(\Phi; s)$ and similarly the volume of
$\mathcal{V}_N(\Phi)$ appears as a special value of
$\widehat{h_N}(\Phi; s)$.
\begin{theorem}
\label{thm:3}
The volume of $\mathcal{U}_N(\Phi)$ is given by
\[
\lambda_{N+1}(\mathcal{U}_N(\Phi)) = 2\widehat{f_N}(\Phi; N+1),
\]
and the volume of $\mathcal{V}_N(\Phi)$ is given by
\[
\lambda_{2N+2}(\mathcal{V}_N(\Phi)) = 2\pi\widehat{h_N}(\Phi; 2N+2).
\]
\end{theorem}

Beyond the computation of the volumes of $\mathcal{U}_N(\Phi)$ and
$\mathcal{V}_N(\Phi)$, any analytic continuation of
$\widehat{f_N}(\Phi; s)$
and $\widehat{h_N}(\Phi; s)$ beyond the range of convergence may yield
further information about the range of values of $\Phi$ which may not
be realizable from other methods.

It should be remarked that the proofs of Theorem~\ref{thm:9} and
Theorem~\ref{thm:3} do not rely on the multiplicativity of $\Phi$.

It is not obvious that, for any choice of $\Phi$, the integrals
defining $\widehat{f_N}(\Phi; s)$ and $\widehat{h_N}(\Phi; s)$ can be
expressed in terms of well-known analytic functions.  As a first step
in this direction, we view the integral defining $\widehat{f_N}$ as a
Lebesgue-Stieltges integral and use integration by parts to write
\begin{equation}
\label{eq:44}
\widehat{f_N}(s) = \left.-\frac{\xi^{-s} f_N(\xi)}{s}\right|_0^{\infty} + \frac{1}{s}
\int_0^{\infty} \xi^{-s} \, d f_N(\xi).
\end{equation}
It follows from Theorem~\ref{thm:9} that $f_N(0) = 0$ and that
$f_N(\xi)$ is dominated by $C \xi^N$ for some constant $C$.
Consequently, the first term in (\ref{eq:44}) is $0$.  From the
definition of $d f_N(\xi)$ we can write

\[
\widehat{f_N}(\Phi; s) = \frac{1}{s} F_N(\Phi; s) \qquad \mbox{where}
\qquad F_N(\Phi; s) = \int_{\R^N} \widetilde{\Phi}(\mathbf{b})^{-s} \,
d\lambda_{N}(\mathbf{b}).
\]
Similarly,
\[
\widehat{h_N}(\Phi; 2s) = \frac{1}{2s} H_N(\Phi; s) \qquad \mbox{where}
\qquad H_N(\Phi; s) = \int_{\C^N} \widetilde{\Phi}(\mathbf{b})^{-2s} \,
d\lambda_{2N}(\mathbf{b}).
\]
Both $F_N(\Phi; s)$ and $H_N(\Phi;s)$ converge to analytic functions
in the region $\Re(s) > N$.  We will call these the real and
complex degree $N$ {\it moment functions} of $\Phi$ (respectively). 

\subsection{Examples of Moment Functions}
Chern and Vaaler's original motivation for computing the moment
functions for Mahler's measure was provided by Theorem~\ref{thm:2}.
Amazingly, their computation revealed that $H_N(\mu; s)$ and $F_N(\mu;
s)$ analytically continued to rational functions of $s$ with poles at
integers and a high order zero at $s=0$.  Moreover, they showed that both
$F_N(\mu; s)$ and $\pi^{-N} H_N(\mu; s)$ have rational coefficients.
\begin{theorem}[(S.-J.~Chern, J.~Vaaler)]
\label{thm:7}
$F_N(\mu; s)$ and $H_N(\mu; s)$ analytically continue to rational
functions of $s$.  In particular,
\begin{equation}
\label{eq:4}
H_N(\mu; s) = \frac{\pi^N}{N!} \prod_{n=1}^N \frac{s}{s - n}.
\end{equation}
If $M$ is the integer part of $(N-1)/2$ then 
\begin{equation}
\label{eq:5}
F_N(\mu; s) = \mathcal{C}_N \prod_{m=0}^M \frac{s}{s - (N - 2m)}
\qquad \mbox{where} \qquad \mathcal{C}_N = 2^N \prod_{m=1}^M
\left(\frac{2m}{2m+1} \right)^{N-2m}.
\end{equation}
\end{theorem}
This surprising result provides additional motivation for determining
the moment functions of other multiplicative distance functions.  And
in fact, the author's original motivation for introducing
multiplicative distance functions and their moment functions was to
create a context in which the surprising rational functions identities
of Chern and Vaaler could be explained.  The next result shows that
much of the structure evident in $F_N(\mu; s)$ and $H_N(\mu; s)$
carries over to the moment functions of the reciprocal Mahler's
measure.
\begin{theorem}
\label{thm:8}
$F_N(\rho; s)$ and $H_N(\rho; s)$ analytically continue to rational
functions of $s$.  In particular,
\begin{equation}
\label{eq:6}
H_N(\rho; s) = 2^N \pi^N \prod_{n=1}^N \frac{s}{s^2 - n^2}.
\end{equation}
If $J$ is the integer part of $(N-1)/2$ then
\begin{equation}
\label{eq:7}
F_N(\rho; s) = v_N \prod_{j=0}^J \frac{s^2}{s^2 - (N-2j)^2}, \qquad
\mbox{where} \qquad v_N = \frac{2^N}{N!} \prod_{n=1}^N
\left(\frac{2n}{2n-1}\right)^{N+1-n}.
\end{equation}
\end{theorem}

A variation of (\ref{eq:6}) was established in \cite{sinclair}. 

The parity (evenness/oddness) of $H_N(\rho; s)$ and $F_N(\rho; s)$
should be mentioned.  This symmetry seems to arise from the fact that
$\rho$ is the pullback of Mahler's measure through the polynomial
$x+1/x$.  The Mellin transform translates the symmetry $x \mapsto 1/x$
to the observed parity in the moment functions.  We may view the
parity of $H_N(\rho; s)$ and $F_N(\rho; s)$ as a kind of functional
equation, and it seems likely that the mechanism which produces this
functional equation will produce functional equations for moment
functions for other multiplicative distance functions formed from the
pullback of Mahler's measure through other rational functions.

As $\mu$ and $\rho$ are the `endpoints' of a `path' of multiplicative
distance functions so are $F_N(\mu; s)$ and $F_N(\rho; s)$ the
`endpoints' of a `path' of moment functions, and similarly for $H_N(\mu;
s)$ and $H_N(\rho; s)$.  Much of the structure present in
Theorem~\ref{thm:7} and Theorem~\ref{thm:8} carries over to the moment
functions of $\mu_t$ for $0 < t < 1$.  By investigating the
qualitative properties of the moment functions of $\mu_t$ we may hope
to learn how the structure of moment functions relates to the
underlying multiplicative distance functions, in particular for those
moment functions which arise as pullbacks of Mahler's measure
through rational functions.

\begin{theorem}
\label{claim:1}
Let $0 < t < 1$.  Then, $F_N(\mu_t; s)$ and $H_N(\mu_t; s)$
analytically continue to rational functions of $s$.  In particular,
\[
H_N(\mu_t; s) = \frac{\pi^N s^N}{N!} \prod_{n=1}^N 
\frac{(1 - t^{2n})s + (1 + t^{2n})n }{s^2 - n^2},
\]
and $F_N(\mu_t; s) \in \Q[t](s)$. Moreover if $J$ is the integer part
$(N-1)/2$ then $F_N(\mu_t; s)$ has simple poles at $\pm N, \pm (N-2),
\ldots, \pm (N-2J)$, a zero of multiplicity $J$ at $s=0$ and $J$ other
real zeros on the negative real axis.
\end{theorem}

\begin{figure}[h!]
\centering
\includegraphics[scale=0.65]{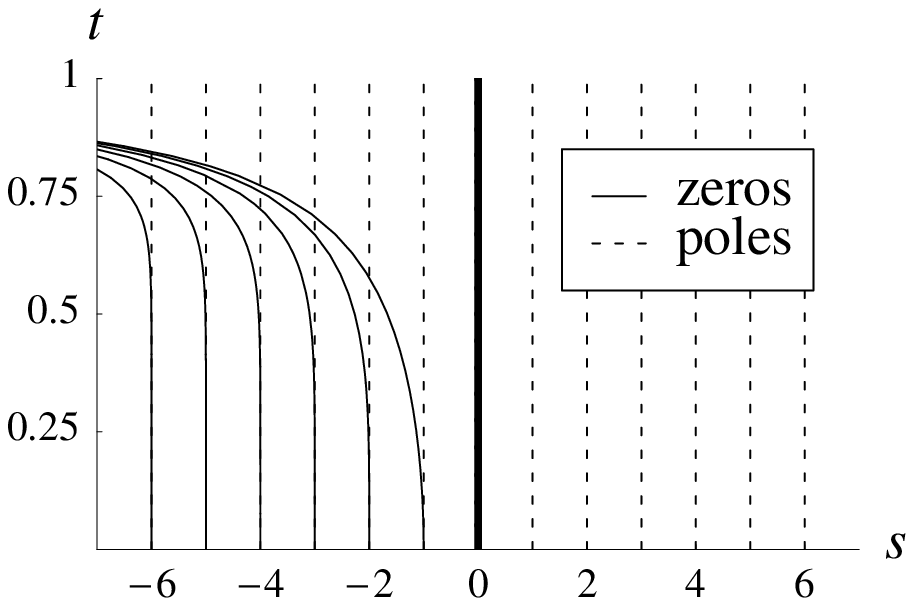}
\includegraphics[scale=0.65]{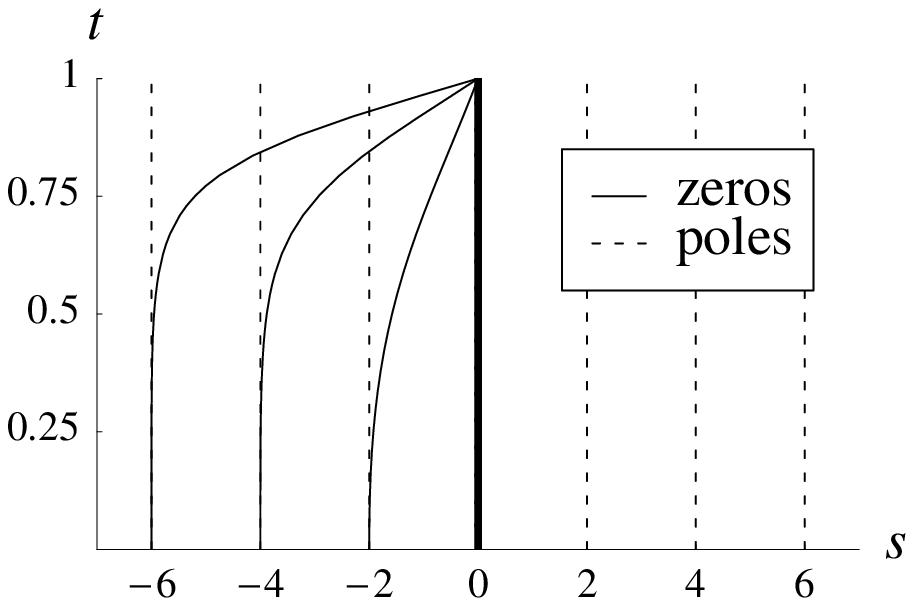}
\begin{caption}{The location of the zeros and poles of $H_6(\mu_t; s)$
  and $F_6(\mu_t; s)$} 
\label{fig:poze}
\end{caption}
\end{figure}

We will leave this theorem unproved since its proof relies on the same
methods we will use to establish Theorem~\ref{thm:8}.  Notice that when
$t=0$ and $t=1$ the formula for $H_N(\mu_t; s)$ coincides with the
formula for $H_N(\mu; s)$ and $H_N(\rho; s)$ respectively.  It should
be remarked that a closed form for $F_N(\mu_t; s)$ can be discovered using
the same method of proof as Theorem~\ref{thm:8}, and this expression
agrees with those for $F_N(\mu; s)$ and $F_N(\rho; s)$ when $t=0$ and
$t=1$.  However, the closed form for $F_N(\mu_t; s)$ is more
complicated than those given for $F_N(\mu; s)$ and $F_N(\rho; s)$ and
in its place we present Figure~\ref{fig:poze}.

In both cases, as $t \rightarrow 0$ the nontrivial zeros
(those not located at $s=0$)
move to cancel the poles located at negative integers.  And, since
$\mu_t$ is the pullback of Mahler's measure through $x + t/x$, this
seems to suggest that the poles at negative integers in $H_N(\mu_t;
s)$ and $F_N(\mu_t;s)$ arise from the $t/x$ factor in $x + t/x$.  
Note that as $t \rightarrow 1$ the nontrivial zeros of
$H_N(\mu_t; s)$ approach $-\infty$ while the nontrivial zeros of
$F_N(\mu_t; s)$ approach 0.  To address this disparity, let
\[
\widetilde{\mathcal{U}_N}(\mu_t) = \{
\mathbf{b} \in \R^N: \widetilde{\mu_t}(\mathbf{b}) = 1 \}
\qquad \mbox{and} \qquad 
\widetilde{\mathcal{V}_N}(\mu_t) = \{
\mathbf{b} \in \C^N: \widetilde{\mu_t}(\mathbf{b}) = 1 \}.
\]
The definitions of $H_N(\mu_t; s)$ and $F_N(\mu_t; s)$ imply that 
\begin{equation}
\label{eq:61}
\lim_{s \rightarrow \infty} F_N(\mu_t; s) =
\lambda_N(\widetilde{\mathcal{U}_N}(\mu_t))
\qquad \mbox{and} \qquad
\lim_{s \rightarrow \infty} H_N(\mu_t; s) =
\lambda_{2N}(\widetilde{\mathcal{V}_N}(\mu_t)).
\end{equation}
From the definition of $\mu_t$ we see that $\mathbf{b}$ is in
$\widetilde{\mathcal{U}_N}(\mu_t)$ exactly when $x^N + \sum_{n=1}^N b_n
x^{N-n}$ has all of its roots in the elliptical region $E_t$.  As $t
\rightarrow 1$, $E_t$ approaches the interval $[-2, 2]$ on the real
axis and the volume of $\widetilde{\mathcal{V}_N}(\mu_t)$ approaches 0
since if $\mathbf{b} \in \widetilde{\mathcal{V}_N}(\rho)$ then in
fact $\mathbf{b} \in \widetilde{\mathcal{U}_N}(\rho)$.  If we
momentarily identify $\C^N$ with $\R^{2N}$ then we see that
$\widetilde{\mathcal{V}_N}(\rho)$ is a subset of codimension $N$ in
$\R^{2N}$.  It is exactly this fact which explains why $H_N(\rho; s)$
has $2N$ poles and only $N$ zeros.  Moreover the fact that as $t
\rightarrow 1$ the non-trivial zeros of $H_N(\mu_t; s)$ tend toward
$-\infty$ verifies our intuition that the volume of
$\widetilde{\mathcal{V}_N}(\mu_t)$ is tends toward 0.  On the other
hand, since $\widetilde{\mathcal{U}_N}(\rho)$ has positive
$\lambda_N$-measure we expect $F_N(\rho; s)$ to have the same number
of zeros and poles, which explains why the non trivial zeros of
$F_N(\mu_t; s)$ do not tend toward $-\infty$ as $t \rightarrow 1$.
The fact that these zeros tend toward $s=0$ seems to support the
hypothesis that the evenness of $F_N(\rho; s)$ stems from the
invariance of $x + 1/x$ under the map $x \mapsto 1/x$.

We remark that explicit formulae for $h_N(\mu_t; \xi)$ and $f_N(\mu_t;
\xi)$ may be recovered from $H_N(\mu_t; s)$ and $F_N(\mu_t; s)$ via
Mellin inversion.  In lieu of explicit formulae we give the following
qualitative corollary Theorem~\ref{claim:1}.  This corollary follows
immediately from Mellin inversion and we will not prove it here.
\begin{corollary}
For each $t \in [0, 1]$, $f_N(\mu_t; \xi)$ and $h_N(\mu_t; \xi)$ are
Laurent polynomials. Moreover $f_N(\mu_t; \xi)$ and $\pi^{-N}
h_N(\mu_t; \xi)$ are in 
$\Q[\xi, \xi^{-1}]$.  
\end{corollary}

\subsection{The Number of  Reciprocal Polynomials in $\Z[x]$ with
  Bounded Degree and Mahler Measure} 
We now turn to an application of the theory of multiplicative distance
functions to Diophantine geometry.  The mechanism by which we may
infer information about the range of $\rho$ on $\Z[x]$ from the range
of $\rho$ on $\R[x]$ stems from the fact that as $T \rightarrow
\infty$ the cardinality of $\Z^{N+1} \cap T \mathcal{U}_N$ is
approximately the volume of $T \mathcal{U}_N$.  Of course we may apply
this principle more generally, but we limit ourselves to the case
of the reciprocal Mahler's measure since reciprocal polynomials hold a
distinguished role in the study of integer polynomials with small
Mahler's measure \cite{smyth}.

A polynomial of degree $M$ is called reciprocal if $f(x) = x^M
f(1/x)$.  Each reciprocal polynomial in $\Z[x]$ corresponds to a
reciprocal Laurent polynomial in $\Z[x+1/x]$.  We denote the set of
reciprocal polynomials in $\Z[x]$ with degree at most $N$ and Mahler's
measure less than or equal to $T$ by $\mathcal{M}_N(T)$.

\begin{theorem}
\label{thm:10}
Let $N$ be a positive integer.  Then, as $T \rightarrow \infty$, the
cardinality of $\mathcal{M}_N(T)$ satisfies the following asymptotic
estimates.  
$$
\# \mathcal{M}_N(T) = \left\{ 
\begin{array}{ll}
\lambda_{J+1}\left(\mathcal{U}_J(\rho)\right) T^{J+1} + O(T^J) &
\mbox{if} \; N=2J, \\
2 \lambda_{J+1}\left(\mathcal{U}_{J}(\rho)\right) T^{J+1} + O(T^{J}) &
\mbox{if} \; N=2J+1.
\end{array}
\right.
$$
where the constant implicit in the $O$-notation is dependent on $J$.
\end{theorem}

\begin{proof}
Suppose that $f$ is a reciprocal polynomial in $\Z[x]$.  The subset of
reciprocal polynomials of $\Z[x]$ is closed under multiplication.  If
$\deg(f)$ is odd then $f(-1) = 0$, and $f(x)/(x+1)$ is a reciprocal
polynomial of even degree.  Furthermore, the multiplicativity of
Mahler's measure implies
$$
\mu(f) = \mu\left(\frac{f(x)}{x+1}\right).
$$
Thus, when studying the range of values Mahler's measure takes on
reciprocal polynomials, it suffices to consider only even degrees.  We
assume that $\deg(f) = 2J$, and let $p(x) = x^{-J} f(x)$.  Clearly $p$
is a reciprocal Laurent polynomial and there exists $g(x) \in \Z[x]$ 
such that $p(x) = g(x + 1/x)$.   It follows that $\mu(f) = \rho(g)$.  

We now turn to $\mathcal{M}_N(T)$.  Notice that $\mathcal{M}_N(T)$
consists of polynomials with both even and odd degrees.  By our
previous remarks, if $N = 2J+1$ is odd, then the set of polynomials in
$\mathcal{M}_N(T)$ with odd degree is in one to one correspondence
with the set
$$
\left\{\mathbf{a} \in \Z^{J}: \rho(\mathbf{a}) \leq T \right\} = T
\mathcal{U}_{J}(\rho) \cap \Z^{J}.
$$
Likewise the set of polynomials in $\mathcal{M}_N(T)$
with even degree is in one to one correspondence $T
\mathcal{U}_{J}(\rho) \cap \Z^{J}$.  

If $N=2J$ is even, the set of polynomials in $\mathcal{M}_N(T)$
with odd degree is in one to one correspondence with the set
$T \mathcal{U}_{J-1}(\rho) \cap \Z^{J-1}$,
while the set of polynomials in $\mathcal{M}_N(T)$
with even degree is in one to one correspondence with the set
$T \mathcal{U}_{J}(\rho) \cap \Z^{J}$.

Then, by well known results, the Lebesgue measure of $T
\mathcal{U}_{J}(\rho)$ gives a good approximation of the number of
integer lattice points contained in $T \mathcal{U}_{J}(\rho)$ when
$T$ is large.  Specifically,
$$
\#\left(T \mathcal{U}_{J}(\rho) \cap \Z^{J}\right) =
\lambda_{J+1}(\mathcal{U}_J(\rho)) T^{J+1} + O(T^J) \quad \mbox{as}
\quad T \rightarrow \infty.
$$
See \cite[Chapter VI, \S2]{lang} or \cite[\S12]{chern-vaaler}
for details.  
\end{proof}

This Theorem is useful since we can explicitly compute
$\lambda_{J+1}(\mathcal{U}_J(\rho))$ using Theorem~\ref{thm:3}.  For example,
\[
\begin{array}{ll}
 \# \mathcal{M}_0(T) = 2 T + O(1), & 
 \# \mathcal{M}_1(T) = 4 T + O(1), \\ & \\
 \# \mathcal{M}_2(T) = \frac{16}{3} T^2 +
  O(T), & 
 \# \mathcal{M}_3(T) = \frac{32}{3} T^2 +
  O(T),  \\ & \\
 \# \mathcal{M}_4(T) = \frac{64}{5} T^3 +
  O(T^2), &
 \# \mathcal{M}_5(T) = \frac{128}{5} T^3 +
  O(T^2), \\& \\
 \# \mathcal{M}_6(T) = \frac{131072}{4725} T^4 +
  O(T^3), & 
 \# \mathcal{M}_7(T) =
  \frac{262144}{4725} T^4 + 
  O(T^3), \\ & \\
 \# \mathcal{M}_8(T) = \frac{655360}{11907} T^5 +
  O(T^4), &
 \# \mathcal{M}_9(T) =
  \frac{1310720}{11907} T^5 + O(T^4), \\ & \\
 \# \mathcal{M}_{10}(T) =
  \frac{2147483648}{21223125} T^6 + O(T^5), & 
 \# \mathcal{M}_{11}(T) =
  \frac{4294967296}{21223125} T^6 + O(T^5).  
\end{array}
\]

\subsection{The Structure of Moment Functions}\label{sec:moment-functions}

The evaluation of $F_N(\Phi; s)$ and $H_N(\Phi; s)$ depends on the
multiplicativity of $\Phi$ as well as the specifics of the root
function $\phi$.  By exploiting the multiplicativity of $\Phi$ we may
express $F_N(\Phi; s)$ and $H_N(\Phi; s)$ in fairly simple terms
dependent only on $\phi$ and $N$ (and of course $s$).  For the
remainder of this section we will view $\Phi$ as fixed.  Many of the
structures introduced in this section are dependent on $\Phi$, but 
this dependence will be suppressed in an effort to simplify the notation.

\begin{sloppypar}

We begin with $H_N(s)$.  For each $s = \sigma + i t$ with $\sigma > N$, let
$\eta_{s}$ be the Borel measure on $\C$ defined by
\[
d\eta_{s}(\gamma) = \phi(\gamma)^{-2 \sigma} \,d\lambda_2(\gamma). 
\]
Next we define a Hermitian form on the Hilbert space $L^2(\eta_s)$ by setting
\[
\la P | Q \ra = \int_{\C} \phi(\gamma)^{-2s} P(\gamma)
\overline{Q(\gamma)} \, d\lambda_2(\gamma) \quad \mbox{for each} \quad
P, Q \in L^2(\eta_s). 
\]
Notice that when $s$ is real this is just the inner product
associated to the norm on $L^2(\eta_s)$.  
It is easy to verify from Theorem \ref{thm:1} that the polynomials $1 ,
\gamma , \gamma^2 , \ldots , \gamma^{N-1}$ are in $L^2(\eta_s)$.  In
fact, any {\it complete set} of $N$ polynomials, that is a set $\{P_1,
P_2, \ldots, P_N\}$  in $\C[\gamma]$ with $\deg P_n = n-1$, is in
$L^2(\eta_s)$.
\end{sloppypar}
\begin{theorem}
\label{thm:4}
Let $\Re(s) > N$, and let $\mathbf{P} = \{P_1, P_2, \ldots,
P_N\}$ be any complete set of monic polynomials in $\C[\gamma]$.  Then,
\[
H_N(\Phi; s) = \det W_{\mathbf{P}}, 
\]
where $W_{\mathbf{P}}$ is the $N \times N$ matrix whose $j,k$ entry is given by
$W_{\mathbf{P}}[j,k] = \la P_j| P_k\ra$.
\end{theorem}

The matrix $W_{\mathbf{P}}$ is known as the Gram matrix of the set
${\mathbf{P}}$ with respect to the Hermitian form $\la \cdot | \cdot
\ra$.  When $s$ is real we may view $\mathbf{P}$ as spanning a
parallelepiped in $L^2(\eta_s)$.  As such $\det H_N(s)$ is the volume
of this parallelepiped.  Moreover, since $H_N(N+1)$ is essentially the
volume of the starbody $\mathcal{V}_N$ can also be regarded as the
volume of a parallelepiped in the Hilbert space $L^2(\eta_{N+1})$.

Perhaps the most useful aspect of Theorem \ref{thm:4} is that it is
independent of the complete family of monic polynomials chosen.  Thus,
a wise choice of ${\mathbf{P}}$---for instance one which is
orthogonal with respect to the Hermitian form---may
make the evaluation of $\det W_{\mathbf{P}}$ easy.  Of course the
coefficients of such orthogonal polynomials will be dependent on $s$.
\begin{corollary}
\label{cor:1}
Let $\Re(s) > N$ and let ${\mathbf{Q}} = \{Q_1, Q_2, \ldots, Q_N\}$ be the
complete family of monic polynomials specified by 
\[
\la Q_j | Q_k \ra = \mathfrak{N}_{s}(Q_k) \, \delta_{kj} \qquad
\mbox{for} \qquad j,k=1, \ldots, N.
\]
Then,
\[
H_N(\Phi; s) = \prod_{n=1}^N \mathfrak{N}_{s}(Q_k).
\]
\end{corollary}
When $s$ is real $\mathfrak{N}_{s}(Q_n)$ is simply the norm squared of $Q_n$
in $L^2(\eta_s)$.

As we shall see the evaluation of $F_N( s)$ is much more
complicated, due in part to the fact that a polynomial in $\R[x]$
may have both real and complex roots.  In spite of this difficulty results
similar to Theorem \ref{thm:4} and Corollary \ref{cor:1} are
available.  These can be stated by replacing the Hermitian form used in
the calculation of $H_N(s)$ with a {\it skew-symmetric} bilinear form
associated to $\Phi$.  The matrix of skew-symmetric bilinear forms formed in
analogy with $W_P$ is antisymmetric and we will replace the
determinant with the Pfaffian---an invariant of antisymmetric
matrices---in order to give a succinct formulation of $F_N(s)$.   

In analogy with the Hermitian form $\la \cdot | \cdot \ra$ we introduce
the skew-symmetric bilinear forms $\la \cdot , \cdot \ra_{\R}$ and $\la
\cdot , \cdot \ra_{\C}$ by
\begin{equation}
\label{eq:14}
\la P, Q \ra_{\R} = \int_{\R^2} \phi(x)^{-s} \phi(y)^{-s} P(x) Q(y)
\sgn(y - x) dx \, dy,
\end{equation}
and
\begin{equation}
\label{eq:15}
\la P, Q \ra_{\C} = -2i \int_{\C} \phi(\beta)^{-s}
\phi(\overline{\beta})^{-s} P(\overline{\beta}) Q(\beta)
\sgn \Im(\beta) \, d\lambda_2(\beta),
\end{equation}
where as before these bilinear forms are implicitly dependent on $s$.
The {\it skew} moniker stems from the fact that $\la Q, P \ra_{\R} = -
\la P, Q \ra_{\R}$ (and similarly for $\la \cdot, \cdot \ra_{\C}$).
When $\Re(s) > N$ it is easily verified that the integrals defining
$\la P, Q\ra_{\R}$ and $\la P, Q \ra_{\C}$ converge when $P$ and $Q$
are polynomials of degree at most $N-1$.  We may create another
skew-symmetric bilinear form by specifying that 
\begin{equation}
\label{eq:16}
\la P, Q \ra = \la P, Q \ra_{\R} + \la P, Q \ra_{\C}.
\end{equation}
Now given any complete family of $N$ monic polynomials $\mathbf{P} =
\{P_1, P_2, \ldots, P_N\} \subseteq \C[\gamma]$, we may create the $N
\times N$ antisymmetric matrix $U_{\mathbf{P}}$ whose $j,k$ entry is
given by $\la P_j, P_k \ra$.  As before the entries of this matrix are
functions of $s$.

An important invariant of even rank antisymmetric matrices is the
Pfaffian.  If $N = 2J $ and $U$ is an $N \times N$ antisymmetric
matrix, then the Pfaffian of $U$ is given by
\begin{equation}
\label{eq:52}
\Pf U = \frac{1}{2^J J!} \sum_{\tau \in S_N} \sgn(\tau) \prod_{j=1}^J
U[\tau(2j-1), \tau(2j)],  
\end{equation}
where $S_N$ is the symmetric group on $\{1,2,\ldots, N\}$.  The
Pfaffian is related to the determinant by the formula $\det U = (\Pf
U)^2$ (see for instance \cite[Appendix: Pfaffians]{MR900587}).  Thus,
psychologically at least, the Pfaffian of $U$ may be thought of as the
signed square root of the determinant of $U$.

One of the major results in this manuscript, and the one we will spend
the most time proving, is that $F_N(s)$ can
be represented as the Pfaffian of $U_{\mathbf{P}}$ for any complete set
$\mathbf{P}$ of $N$ monic polynomials in $\R[\gamma]$.  However,
before this claim can be made it is necessary to adjust our
definitions for the case when $N$ is odd.  

\begin{theorem}
\label{thm:5}
\sloppy
Let $\Re(s) > N$ and let $J$ be the integer part of $(N+1)/2$.  If
$\mathbf{P} = \{P_1(\gamma), P_2(\gamma), \ldots, P_N(\gamma)\}$ is
any complete set of monic polynomials in $\C[\gamma]$ then
\[
F_N(\Phi; s) = \Pf U_{\mathbf{P}},
\]
where $U_{\mathbf{P}}$ is the $2J \times 2J$ antisymmetric matrix whose $j,k$
entry is given by 
\begin{equation}
\label{eq:1}
U_{\mathbf{P}}[j,k] = \left\{
\begin{array}{ll}
\la P_j, P_k \ra & \quad \mbox{if } j, k \leq N, \\ 
{\displaystyle 
\sgn(k-j) \int_{\R} \phi(x)^{-s}\, P_{\min\{j,k\}}(x) \, dx
} & \quad \mbox{otherwise}. \\
\end{array}
\right.
\end{equation}
\end{theorem}
Notice that when $N$ is even then the first condition in
equation~(\ref{eq:1}) always holds.

As is the case with $H_N(s)$ a smart choice of $\mathbf{P}$
yields a simple product formulation for $F_N(s)$.  Specifically,
when $N$ is even we may use a complete family of monic polynomials
which are {\it skew-orthogonal}.  
\begin{corollary}
\label{cor:2}
Suppose that $N = 2J$, $\Re(s) > N$ and let $\mathbf{Q} = \{Q_1, Q_2,
\ldots, Q_N\}$ be any complete family of monic polynomials specified
by 
\[
\la Q_{2k-1}, Q_{2j} \ra = -\la Q_{2j}, Q_{2k - 1} \ra = \delta_{kj}
\mathfrak{M}_s(Q_j) \quad \mbox{and} \quad \la Q_{2j}, Q_{2k} \ra = \la
Q_{2j-1}, Q_{2k-1} \ra = 0,
\]for $j,k=1, \ldots, J$.  Then,
\[
F_N(\Phi; s) = \prod_{j=1}^J \mathfrak{M}_s(Q_j).
\]
\end{corollary}
The quantities $\mf{M}_s(Q_j)$ are referred to as the normalization(s)
of $\mathbf{Q}$.  

In the special case of multiplicative distance functions
whose root functions satisfy certain symmetries we may write $F_N(\Phi;
s)$
as a determinant. 
\begin{corollary}
\label{cor:3}
Suppose that $\Re(s) > N$
and let $J$ be the integer part of $(N+1)/2$.   Furthermore suppose
that $\mathbf{P}$ is a
complete family of monic polynomials in $\R[x]$ such that $P_n$ is
even when $n-1$ is even, and $P_n$ is odd when $n-1$ is odd.  If the root function of $\Phi$ satisfies $\phi(-\beta) =
\phi(\beta)$ for every $\beta \in \C$ then,
\[
F_N(\Phi; s) = \det A_{\mathbf{P}}
\]
where $A_{\mathbf{P}}$ is the $J \times J$ matrix whose $j,k$ entry is
given by
\[
A_{\mathbf{P}}[j,k] = U_{\mathbf{P}}[2j-1, 2k].
\]
\end{corollary}

\section{The Proof of Theorem~\ref{thm:1}}
Since $\Phi$ is continuous, non-negative and positive definite, we
find that $\phi$ is continuous and $\phi(\alpha) > 0$ for each $\alpha
\in \C$.  The asymptotic properties of root functions are derived from
the continuity of multiplicative distance functions.  To see this,
let $a$ and $b$ be nonzero complex numbers.  By homogeneity,
$$
\Phi(a x - b) = |a| \, \Phi\left(x - \frac{b}{a}\right) =
|a|\, \phi\left(\frac{b}{a}\right) 
$$ 
By continuity ${\displaystyle \lim_{|a| \rightarrow 0} \Phi(a x - b)
= \Phi(-b) = |b|}$, and thus
$$
\lim_{|a| \rightarrow 0}|a| \,\phi\left(\frac{b}{a}\right) = |b|.
$$
Setting $\gamma = b/a$ we see that $\phi(\gamma) \sim |\gamma|$ as
$|\gamma| \rightarrow \infty$.

The other direction is more complicated.  Suppose that $\psi: \C \rightarrow
(0,\infty)$ is a continuous function such that $\psi \sim |\gamma|$ as
$|\gamma| \rightarrow \infty$.  We will use a modification of
Mahler's original proof that $\mu$ is continuous to prove the
continuity of the function
\[
\Psi: a \prod_{n=1}^N (x - \gamma_n) \mapsto |a| \prod_{n=1}^N \psi(\gamma_n).
\]
Certainly $\Psi$ satisfies the other axioms of multiplicative distance
functions.  

In fact we will prove that $\Psi$ is continuous with respect to the
stronger topology induced by uniform convergence on compact subsets of $\C$.
Suppose that $\{f_k(x)\}$ is a sequence of polynomials in $\C[x]$ such that 
$$
f_k(x) = a_{k N_k} \prod_{n=1}^{N_k} (x - \gamma_{kn}) \quad
\mbox{for} \quad k > 0
$$
and
$$
\lim_{k \rightarrow \infty} f_k(x) = f(x) = a \prod_{n=1}^N (x - \gamma_n),
$$
uniformly on compact subsets of $\C$.  We will show that 
$$
\lim_{k \rightarrow \infty} \Psi(f_k) = \Psi(f).
$$
By an easy corollary to Hurwitz's Root Theorem (see for instance
\cite{saks-zygmund}) we may reorder the roots of each $f_k(x)$ so that 
$$
\lim_{k \rightarrow \infty} \gamma_{kn} = \gamma_n  \quad \mbox{for}
\quad n=1, \ldots, N.
$$
For each $k>0$, define the polynomials $g_k(x)$ and $h_k(x)$ by,
$$ 
g_k(x) = a \prod_{n=1}^N (x - \gamma_{kn}) \quad \mbox{and} \quad
h_k(x) = \frac{a_{kN_k}}{a} \prod_{n=N+1}^{N_k} (x - \gamma_{kn}),
$$ 
and notice that $g_k(x) \rightarrow f(x)$ as $k \rightarrow \infty$.
Furthermore, since $g_k(x)$ is of degree $N$ for all $k$, it follows
that this convergence is uniform on compact subsets of $\C$. Now, 
\begin{eqnarray*}
\left| g_k(x) \right| \left| h_k(x) - 1 \right| &=& \left| f_k(x) -
g_k(x) \right| \\
&\leq& \left| f_k(x) - f(x) \right| + \left| f(x) - g_k(x) \right|,
\end{eqnarray*}
from which it follows that $\{h_k(x)\}$ converges to the constant
function $1$ uniformly on compact subsets of $\C \setminus \{\gamma_1,
\ldots, \gamma_N\}$.  

The constant coefficient of $h_k(x)$ is given by
$$
\frac{a_{k N_k}}{a}(-1)^{N_k - N} \prod_{n=N+1}^{N_k} \gamma_{kn},
$$
and thus, by choosing a point $x \in \C \setminus\{\gamma_1,
\ldots, \gamma_N\}$ and using the fact that $h_k(x) \rightarrow 1$
we have,
\begin{equation}
\label{eq:45}
\lim_{k \rightarrow \infty} \left\{ \frac{\left| a_{N_k} \right|}{\left|
  a \right|} \prod_{n=N+1}^{N_k} \left| \gamma_{kn}
\right| \right\}= 1.
\end{equation}
This is the key fact needed to prove the theorem.
\begin{eqnarray*}
\lim_{k \rightarrow \infty} \Psi(f_k) &=&
\lim_{k \rightarrow \infty} \left( |a_{kN_k}| \prod_{n = 1}^{N_k}
  \psi(\gamma_{kn}) \right) \\
&=& \lim_{k \rightarrow \infty}\left(
\left\{
\frac{|a_{kN_k}|}{|a|} \prod_{n = N+1}^{N_k} \psi(\gamma_{kn})
\right\}
\left\{
|a| \prod_{n=1}^N \psi(\gamma_{kn})
\right\}
\right) \\
&=&
\lim_{k \rightarrow \infty}\left(
\left\{
\frac{|a_{kN_k}|}{|a|} \prod_{n = N+1}^{N_k} |\gamma_{kn}|
\right\}
\left\{
|a| \prod_{n=1}^N \psi(\gamma_{kn})
\right\}
\right).
\end{eqnarray*}
Where the last equation is a consequence of the fact that $|\gamma_{kn}|
\rightarrow \infty$ as $k \rightarrow \infty$, and $\psi(\gamma) \sim
|\gamma|$.  From (\ref{eq:45}) it follows that,
\begin{eqnarray*}
\lim_{k \rightarrow \infty} \Psi(f_k) = \lim_{k \rightarrow \infty}
\left\{ |a| \prod_{n=1}^N \psi(\gamma_{kn}) \right\} 
= |a| \prod_{n=1}^N \psi(\gamma_n) = \Psi(f),
\end{eqnarray*}
where the second equality follows from the continuity of $\psi$ and
the fact that ${\displaystyle \lim_{k \rightarrow \infty} \gamma_{kn}
  = \gamma_n}$.

\section{The Proofs of Theorem~\ref{thm:9} and Theorem~\ref{thm:3}}

We will prove Theorem~\ref{thm:9} for the real case and leave the
complex case to the reader.  We will view $\Phi$ as
fixed and suppress any notational dependence on $\Phi$.

\begin{proof*}[of Theorem~\ref{thm:9}]
Let $B$ be the $N+1$ dimensional unit ball
centered at the origin.  Then, since $\mathcal{U}_N$ is bounded
we can find a positive constant $\eta$ so that
\[
\mathcal{U}_N \subset \eta B \quad \mbox{and thus}
\quad T \mathcal{U}_N \subset T \eta B.
\]
Let $\mathcal{A}_{1/T} = \{(\mathbf{b}, 1/T): \mathbf{b} \in
\R^{N}\}$. For instance, $\mathcal{A}_1$ is the 
hyperplane of coefficient vectors of monic polynomials of degree $N$.
It follows that
\begin{equation}
\label{eq:43}
(\mathcal{A}_1 \cap T \mathcal{U}_N) \subset (\mathcal{A}_1 \cap T \eta \;
B). 
\end{equation}
Depending on the value of $T$, the set $(\mathcal{A}_1 \cap T \eta \;
B)$ is either empty or an $N$-dimensional ball. It follows from
(\ref{eq:43}) that 
\[
f_N(T) \leq \lambda_{N}(\mathcal{A}_1 \cap T \eta B),
\]
and there exists an $\epsilon > 0$ such that if $T < \epsilon$ then $f_N(T) = 0$.

Clearly the set of polynomials with leading coefficient $1/T$ and
distance 1 is given by \mbox{$\mathcal{A}_{1/T} \cap \mathcal{U}_N$}.  Notice that $\mathcal{A}_{1/T} = (1/T) \mathcal{A}_1$.
Thus we find that
\[
\mathcal{A}_{1/T} \cap \mathcal{U}_N = \frac{1}{T}(\mathcal{A}_1 \cap T\mathcal{U}_N).
\]
It is easy to see that $(\mathcal{A}_{1/T} \cap \mathcal{U}_N) \rightarrow
\mathcal{U}_{N-1}$ as $T \rightarrow \infty$.  Thus
\multbox
\begin{eqnarray*}
\lambda_{2N}(\mathcal{U}_{N-1}) &=& \lim_{T \rightarrow
  \infty}\lambda_{N}(\mathcal{A}_{1/T} \cap \mathcal{U}_N) 
\\ &=& \lim_{T 
  \rightarrow \infty} \lambda_{N}\left(\frac{1}{T}(\mathcal{A}_1 \cap
  T\mathcal{U}_N)\right) = \lim_{T \rightarrow \infty}
\frac{f_N(T)}{T^{N}}. 
\end{eqnarray*}
\emultbox
\end{proof*}

Again we will prove Theorem~\ref{thm:3} in the real case and leave the
complex case to the reader.

\begin{proof*}[of Theorem~\ref{thm:3}]
Using the notation from the proof of Theorem~\ref{thm:9}, 
the volume of $\mathcal{U}_N$ is given by
\begin{equation*}
\lambda_{N+1}(\mathcal{U}_N) = \int_{\R} \lambda_{N}(
  \mathcal{A}_T \cap \mathcal{U}_N) \, dT.
\end{equation*}
By the absolute homogeneity of $\Phi$ we see 
\[
\lambda_{N}(\mathcal{A}_T \cap \mathcal{U}_N) = |T|^N
\lambda_N(\mathcal{A}_1 \cap |T|^{-1} \mathcal{U}_N) = |T|^N f_N(|T|^{-1}).
\]
And thus 
\[
\lambda_{N+1}(\mathcal{U}_N)  = \int_{\R} |T|^{N} f_N(|T|^{-1})  \, d T
= 2 \int_0^{\infty} T^N f_N(T^{-1}) \, dT.
\]
Finally, by setting $\xi = T^{-1}$ we find
\[
\singlebox
\lambda_{N+1}(\mathcal{U}_N) = 2 \int_0^{\infty} \xi^{-N-1}\,
f_N(\xi) \, d\xi = 2 \widehat{f_N}(N+1). 
\esinglebox
\]
\end{proof*}

\section{The Proof of Theorem \ref{thm:7}}

As remarked previously Theorem \ref{thm:7} was originally proved by
S-J. Chern and J. Vaaler in \cite{chern-vaaler}.  Their technique for
evaluating $F_N(\mu; s)$ involved a number of rational function
identities which were specialized to $\mu$.  In this section we will
present a different proof which relies on Theorem \ref{thm:5} (or
rather Corollary \ref{cor:3}).  

It is worth remarking that Theorem~\ref{thm:5} (and its corollaries)
reduce the determination of $H_N(\Phi; s)$ and $F_N(\Phi; s)$ to the
(not necessarily trivial) calculation a number of Hermitian forms and
skew-symmetric bilinear forms dependent on $\Phi$, and then the
computation of a determinant and a Pfaffian.  When $\Phi = \mu$ it is
convenient to use the family of monic polynomials $\mathbf{P} = \{1,
\gamma, \gamma^2, \ldots, \gamma^{N-1}\}$.  For $H_N(\mu;
s)$ the integrals defining the Hermitian forms of pairs of elements of
$\mathbf{P}$ are elementary, and moreover $\mathbf{P}$ is orthogonal
with respect to the Hermitian form.  Thus, $H_N(\mu; s)$ is the
determinant of a diagonal matrix with entries that are easily
computed.  The details of this computation are left to the reader (or
can be found in \cite{chern-vaaler}).

The integrals defining the skew-symmetric bilinear forms used in the
computation of $F_N(\mu; s)$ are slightly more complicated, but still
elementary.  And since $\mu$ satisfies the conditions of
Corollary~\ref{cor:3}, $F_N(\mu; s)$ is given by the determinant of a
matrix whose entries are given by these skew-symmetric bilinear 
forms.  This matrix is more complicated than the matrix which appears
in the formulation of $H_N(\mu; s)$ but nonetheless its determinant
can be computed.
\begin{lemma}
\label{lemma:1}
Let $\Phi = \mu$ and set $J$ be the integer part of $(N+1)/2$. Define
$A_{\mathbf{P}}$ to be the $J \times J$ matrix defined as in Corollary \ref{cor:3}.  Then,
\begin{equation}
\label{eq:17}
A_{\mathbf{P}}[j,k] = \left\{
\begin{array}{ll}
{\displaystyle
 \frac{1}{2k-2j+1}\left(\frac{4}{2j-1}\right)  \left(\frac{s}{s-2k}\right)}
 & \quad \mbox{if} \quad {\displaystyle k < \frac{N+1}{2}}, \\ &  \\
{\displaystyle 
\frac{2}{2j-1}\left(\frac{s}{s-2j+1}\right)
} & \quad \mbox{if} \quad {\displaystyle k = \frac{N+1}{2}}. \\
\end{array}
\right.
\end{equation}
\end{lemma}

We will defer the proof of this lemma to see how we may use it to
derive the formulation for $F_N(\mu; s)$ given in Theorem~\ref{thm:7}.
Since the second condition defining $A_{\mathbf{P}}[j,k]$ is only
realized when $N$ is odd, it is sensible to divide the determination
of $F_N(\mu; s)$ into cases depending on whether $N$ is even or odd.

\subsection{The Even $N$ Case}
When $N = 2J$, $A_{\mathbf{P}}$ is defined only by the first condition
in (\ref{eq:17}).  We have written this suggestively to indicate terms
which depend only on the rows or columns of $A_{\mathbf{P}}$.  It
follows that
\[
F_N(\mu; s) = \det A_{\mathbf{P}}[j,k] = \det B \cdot 2^N \left\{ \prod_{j=1}^J
\frac{s}{s-2j}\left(\frac{1}{2j-1}\right) \right\},
\]
where $B$ is the $J \times J$ matrix given by $B[j,k] = 1/(2k-2j+1)$.
The matrix $B$ is a Cauchy matrix, and using the well-known formula
for the determinant of a Cauchy matrix, 
\begin{equation}
\label{det B}
\det B = (-1)^{J \choose 2} \left. \left\{\prod_{1 \leq j < k \leq
  J}\!\!\! (2k - 2j)^2 \right\} \right/ \left\{
\prod_{j=1}^J \prod_{k=1}^J (2k - 2j + 1)
\right\},
\end{equation}
the denominator of which is
\begin{eqnarray*}
\prod_{j=1}^J \prod_{k=1}^J (2k - 2j + 1) &=& (-1)^{J \choose 2}
\left\{ \prod_{1 \leq j < k \leq  J}\!\!\! (2(k - j) + 1) (2(k - j) -
1)  \right\} \\
&=& (-1)^{J \choose 2} \left\{ \prod_{1 \leq j < k \leq  J}\!\!\! (2(k
- j) + 1)^2 \right\} \left\{\prod_{j=1}^J (2(J-j) + 1)\right\}.
\end{eqnarray*}
Substituting this into (\ref{det B}) we find
\begin{equation}
\label{eq:19}
\left\{\prod_{j=1}^J \frac{1}{2j-1} \right\} \cdot \det B = 
\left\{ \prod_{j=1}^{J-1} \left(\frac{2 j}{2j + 1}\right)^{2J-2j} \right\}.
\end{equation}
And thus,
$$
F_N(\mu; s) = 2^N \left\{\prod_{j=1}^{J-1} \left(\frac{2 j}{2j +
      1}\right)^{N-2j} \right\} \left\{ \prod_{j=1}^J \frac{s}{s - 2j}
\right\},
$$
which after reindexing yields the formula for $F_N(\mu; s)$ given
in Theorem \ref{thm:7}.

\subsection{The Odd $N$ Case}
When $N$ is odd we have $J = (N+1)/2$.  Looking at (\ref{eq:17}) we
may factor out terms dependent only on the rows or columns of $A_{\mathbf{P}}$ to write
\begin{equation}
\label{eq:18}
F_N(\mu; s) = \det A_{\mathbf{P}} = 2^N s^J  \left\{
  \prod_{j=1}^J \frac{1}{2j-1}\right\} \left\{ \prod_{k=1}^{J-1}
  \frac{1}{s - 2k} \right\} \cdot \det B',
\end{equation}
where $B'$ is the $J \times J$ matrix given by
\[
B'[j, k] = \left\{
\begin{array}{lc}
{\displaystyle \frac{1}{2k - 2j + 1}} & \quad \mbox{if} \quad k < J, \\
& \\
{\displaystyle \frac{1}{s - 2j + 1}} & \quad \mbox{if} \quad k = J.
\end{array}
\right.
\]
The determinant of $B'$ is clearly a rational function of $s$ which we
will denote by $b(s)$.  Moreover as $s \rightarrow \infty$ we must
have $b(s) \rightarrow 0$, from which it follows that $b(s)$ has fewer
zeros than poles.  It is clear from the definition of $B'$ that $b(s)$
has $J$ simple poles located at the positive odd integers not
exceeding $N$.  It is also easy to see that $b(s)$ has $J-1$ zeros
located at the positive even integers not exceeding $N$.  Thus there
exists a constant $\kappa$ such that
\[
b(s) = \kappa \left\{ \prod_{j=1}^J \frac{1}{s - 2j + 1} \right\}
\left\{  \prod_{k=1}^{J-1} (s - 2k) \right\}
\]
Notice that the zeros of $b(s)$ exactly cancel the poles at even
integers which appear in (\ref{eq:18}).  That is,
\begin{equation}
\label{eq:20}
F_N(\mu; s) = \kappa \cdot 2^N  \left\{ \prod_{j=1}^J
\frac{s}{s - 2j + 1} \left(\frac{1}{2j-1}\right) \right\}.
\end{equation}In order to determine the value of $\kappa$ we must determine $b(s)$ at
another value of $s$, the obvious choice being $s = 2J$.  In this
situation $b(2J)$ is simply the determinant of the $J \times J$ matrix
whose $j,k$ entry is given by $1/(2k - 2j + 1)$.  That is, $b(2J) =
\det B$, the same Cauchy determinant that appeared in the even $N$ case.  Thus we have
\[
b(2J) =  \kappa \left\{ \prod_{j=1}^{J-1} \frac{2J - 2j}{2J - 2j + 1}
\right\} = \det B,
\]
and by (\ref{eq:19}) 
\[
\kappa = \left\{ \prod_{j=1}^J (2j - 1) \right\} \left\{
  \prod_{j=1}^{J-1} \left( \frac{2j }{2j+1}
  \right)^{2J - 2j - 1} \right\}.
\]
Substituting this into (\ref{eq:20}) we find
\[
F_N(\mu; s) = 2^N   \left\{\prod_{j=1}^{J-1} \left( \frac{2j }{2j+1}
  \right)^{N - 2j} \right\} \left\{ \prod_{j=1}^J \frac{s}{s- 2j + 1} \right\},
\]
which after reindexing yields the formulation for $F_N(\mu; s)$ given
in Theorem~\ref{thm:7}.

\subsection{The Proof of Lemma~\ref{lemma:1}}

We compute the entries of the matrix $U_{\mathbf{P}}[j,k] = \la
\gamma^{j-1}, \gamma^{k-1} \ra$ under
the conditions that $j$ is odd and $k$ is even.  
  The root function of $\mu$ is $\phi(\gamma) = \max\{1, |\gamma|\}$,
  and hence
\begin{eqnarray*}
\la \gamma^{j-1}, \gamma^{k-1} \ra_{\R} &=&
  \int_{-\infty}^{\infty} 
\int_{-\infty}^{\infty} \max\{1, |x| \}^{-s} \max\{1,
|y|\}^{-s} x^{j-1} y^{k-1} \sgn(y - x) \, dx \, dy   \\ 
&=& 2
\int_{-\infty}^{\infty} \int_{-\infty}^{y} \max\{1, |x|
\}^{-s} \max\{1, |y|\}^{-s} x^{j-1} y^{k-1} \, dx \, dy.
\end{eqnarray*}
But this integral is elementary, since we may divide the domain of
integration into regions according to where $\max\{1, |x|\}$ and
$\max\{1, |y|\}$ are identically one.  The integrals converge when
$\Re(s) > j + k$.  Putting the result into partial fractions form (as
a function of $s$) we find,
\begin{equation}
\label{eq:8}
\la \gamma^{j-1},\gamma^{k-1} \ra_{\R} = \underbrace{\frac{2}{2s - j - k}
\left( \frac{4}{j-k} \right)}_{\mding{182}} + \frac{4}{j(j+k)} + \frac{2}{s -
  k}\left(\frac{2k}{j(k-j)} \right).
\end{equation}

Now,
\begin{eqnarray*}
\la \gamma^{j-1}, \gamma^{k-1} \ra_{\C} &=& -2i \int_{\C} \max\{1,
|\beta|\}^{-2s} 
(\overline{\beta})^{j-1} \beta^{k-1} \sgn \Im(\beta) \,
d\lambda_2(\beta)  \\
&=& -2 i \int_0^{\infty} \max\{1, r\}^{-2s} r^{j+k-1} \, dr
\times \left\{ \int_0^{\pi} - \int_{\pi}^{2 \pi} \right\} e^{(k - j)i \theta}
\, d\theta \\
&=& -4 i \int_0^{\infty} \max\{1, r\}^{-2s} r^{j+k-1} \, dr
\times \int_0^{\pi} e^{(k - j)i \theta}\, d\theta.
\end{eqnarray*}
The integrals in this expression are elementary, and when $\Re(s)
> j+ k$, we find
\[
\la \gamma^{j-1}, \gamma^{k-1} \ra_{\C} = -\frac{2}{2s - j -
  k}\left(\frac{4}{j-k}\right) + \frac{8}{(k-j)(k+j)}
\]
Notice that the first term in $\la \gamma^{j-1}, \gamma^{k-1}
\ra_{\C}$ exactly cancels \ding{182} in (\ref{eq:8}).  That is,
\begin{align*}
A_{\mathbf{P}}[j,k] = \la \gamma^{j-1},\gamma^{k-1} \ra &= \frac{4}{j(j+k)} 
+ \frac{2}{s - k}\left(\frac{2k}{j(k-j)} \right) +
\frac{8}{(k-j)(k+j)} \nonumber \\ &= \displaystyle \frac{4}{j(j-k)}
\left(\frac{s}{k-s}\right). 
\end{align*}
When $k \leq N$, the entries of $A_{\mathbf{P}}[j,k]$ are given by
$U_{\mathbf{P}}[2j-1, 2k]$.  When $k = N + 1$ (which can only
occur when $J$ is odd) we have 
\[
A_{\mathbf{P}}[j, N+1] = \int_{-\infty}^{\infty} \max\{1,
|x|\}^{-s} x^{2j-2} \, dx = \frac{2 s}{(2j-1)(s-2j+1)},
\]
which finishes the proof of the lemma.  

\section{The Proof of Theorem \ref{thm:8}}
Our strategy is the same as in the proof of Theorem~\ref{thm:7}:  First
compute the Hermitian and skew-symmetric bilinear forms for a
complete family of polynomials and then compute the determinant and
Pfaffian of the appropriate matrices whose entries are these bilinear
forms.  To compute $F_N(\rho; s)$  we will use Corollary~\ref{cor:3} with $\mathbf{P}=\{1, \gamma,
\gamma^2, \ldots, \gamma^{N-1}\}$.   The formulation of $H_N(\rho; s)$
given in Theorem~\ref{thm:8} will not be presented here since it (or
rather a minor variation of it) is the subject of \cite{sinclair}.

\begin{lemma}
\label{lemma:2}
Let $\Phi = \rho$ and let $J$ be the integer part of $(N+1)/2$.
Define $A_{\mathbf{P}}$ to be the $J \times J$ matrix
defined as in Corollary~\ref{cor:3}.  Then, if $k < (N+1)/2$,
\begin{equation}
\label{eq:21}
A_{\mathbf{P}}[j,k] = 
  \sum_{n=1}^{J} \left[ {2k-1 \atop k - n}
    \right] \left(\frac{16 s^2}{s^2 - (2n)^2} \right) 
  \sum_{m=1}^{J} \left[ {2j-2 \atop j - m}
    \right] \frac{2n}{2 m - 1}
    \left(\frac{1}{(2n)^2 - (2 m-1)^2} \right),
\end{equation}
and,
\begin{equation}
\label{eq:37}
A_{\mathbf{P}}[j,{\textstyle \frac{N+1}{2}}] = \sum_{n=1}^J \left[{N \atop \frac{N+1}{2}-n} \right] \frac{s^2}{2^{N-2}}
\sum_{m=1}^J \left[{2j-2 \atop j-m} \right] \frac{2n}{2m-1}
\left(\frac{1}{s^2 - (2m-1)^2} \right),
\end{equation}
where
\begin{equation}
\label{eq:22}
\left[ M \atop m \right] = 
{M \choose m} - {M \choose m-1}.
\end{equation}
\end{lemma}

Lemma~\ref{lemma:2} is proved by brute force, and in fact much of the
proof involves massaging the entries of $A_{\mathbf{P}}$ into the form
given in the statement of the lemma.  This form is not the
most natural but will be useful for our purposes as we ultimately need to
take the determinant of $A_{\mathbf{P}}$.  We defer the proof of
Lemma~\ref{lemma:2} in order to see how it may be used to compute the
formulation of $F_N(\rho; s)$ given in the statement of
Theorem~\ref{thm:8}.

We define the $J \times J$ matrices $B, C$ and $D$ by
\[
B[m,n] = \left\{
\begin{array}{ll}
{\displaystyle \frac{2n}{2m-1} \left(\frac{1}{(2n)^2 - (2m-1)^2}\right)} & \mbox{if}
\quad n < \frac{N+1}{2}, \\ & \\
{\displaystyle \frac{2n}{2m-1} \left(\frac{1}{s^2 - (2m-1)^2} \right)}
  & \mbox{if} 
\quad n = \frac{N+1}{2},
\end{array}
\right. 
\]
\[
C[j,m] = \left[{2j-2 \atop j-m} \right] \qquad \mbox{and} \qquad D[k,n] = \left\{
\begin{array}{ll}
{\displaystyle 
\left[{2k-1 \atop k-n} \right] \left(\frac{16 s^2}{s^2 - (2n)^2} \right)}
& \mbox{if} \quad k < \frac{N+1}{2}, \\ & \\
\displaystyle{
\left[{2k-1 \atop k-n} \right] \frac{s^2}{4^k}} & \mbox{if} \quad k = \frac{N+1}{2}.
\end{array}
\right.
\]Thus $A_{\mathbf{P}} = DCB$ and $F_N(\rho; s) = \det A_{\mathbf{P}}
= \det B \cdot \det C \cdot \det D$.  This is convenient since, the
matrices $C$ and $D$ are triangular (since for instance if $m > j$
then $\left[ {2j-2 \atop j-m}\right] = 0$).

When $N$ is even we have $J < (N+1)/2$ and hence the first conditions
defining $B$ and $D$ hold.  Thus $\det B$ is simply a rational number,
and thus computing the diagonal entries of $C$ and $D$ we see that
there is some rational number $v_N$ so that
\begin{equation}
\label{eq:38}
F_N(\rho; s) = v_N \prod_{j=1}^J \frac{s^2}{s^2 - (2j)^2} =
v_N \prod_{j=0}^{J} \frac{s^2}{s^2 - (N-2j)^2}
\end{equation}
When $N$ is odd, then the determinant of $B$ is a rational function of $s$,
which we will denote $b(s)$.  From the definition of $B$ it is easily
seen that $b(s)$ is an even rational function with simple poles at the
integers $\pm 1, \pm 3, \ldots, \pm N$.  Moreover when $s = \pm 2, \pm
4, \ldots, \pm (N-1)$ the matrix $B$ is singular and hence $b(s) = 0$
for these values of $s$.  Also from the definition of $B$ it is seen
that $b(s) \rightarrow 0$ as $s \rightarrow \infty$.  We conclude that
there $b(s)$ has fewer zeros than poles.  We have identified all the
poles of $b(s)$, and since $b(s)$ is even we have also identified the
complete list of zeros of $b(s)$.  Notice that the zeros of $b(s)$
exactly cancel the poles which arise from the diagonal entries of
$D$.  Putting these observations together we find that there exists
some rational number (which we also denote $v_N$) such that
\begin{equation}
\label{eq:39}
F_N(\rho; s) = v_N \prod_{j=}^J \frac{s^2}{s^2 - (2j-1)^2} = v_N
\prod_{j=0}^J \frac{s^2}{s^2 - (N-2j)^2}.
\end{equation}
We may find the value of $v_N$ by explicitly computing the determinant
of $B$, $C$ and $D$, and noting that $\det B$ is a Cauchy
determinant.  However, by casting the constant $v_N$ in another
context we may find its value in the literature.

Equation~(\ref{eq:61}) implies that $v_N =
\lambda_N(\widetilde{\mathcal{U}_N}(\rho))$.  That is $v_N$ is the
volume of the set of $\mathbf{b} \in \R^N$ with
$\widetilde{\rho}(\mathbf{b}) = 1$.  This observation is useful since
the volume of $\widetilde{\mathcal{U}_N}(\rho)$ has been computed by 
S.~DiPippo and E.~Howe in
\cite[Proposition~2.2.1]{dipippo-howe}. DiPippo and Howe show that
\begin{equation}
\label{eq:33}
\lambda_N(\widetilde{\mathcal{U}_N}(\rho)) = \frac{2^N}{N!} \prod_{n=1}^N
\left(\frac{2n}{2n-1} \right)^{N+1-n}.
\end{equation}In fact, DiPippo and Howe report this number as the volume of the set of monic
coefficient vectors of polynomials of degree $N$ with real
coefficients and all roots on the unit circle, however in the course
of their computation they show that this volume exactly equals the
volume of $\lambda_N(\widetilde{\mathcal{U}_N}(\rho))$.  Putting
(\ref{eq:33}) together with (\ref{eq:38}) and (\ref{eq:39}) we arrive
at the formulation of $F_N(\rho; s)$ given in the statement of
Theorem~\ref{thm:8}. 

\subsection{The Proof of Lemma~\ref{lemma:2}}
Throughout this section we set $\Phi =
\rho$.

Before proving Lemma \ref{lemma:2} a few results about the
binomial-like coefficients are in order.  First let us see how these
coefficients come about.  The entries of $A_{\mathbf{P}}$ are defined
by integrals in which factors like $\phi(\gamma)^{-s} \gamma^{j-1}$
occur.  In order to evaluate these integrals it is convenient to use
the change of variables $\gamma \mapsto \gamma + 1/\gamma$ since
$\phi(\gamma + 1/\gamma) = \max\{|\gamma|, |\gamma|^{-1} \}$.  We are
left with integrands of the form 
\[
\max\{|\gamma|, |\gamma|^{-1}\}^{-s} \left( \gamma + \frac{1}{\gamma}
\right)^{j-1} \left|1 - \frac{1}{\gamma^2}  \right|.
\]
This is beneficial since we may use the
binomial theorem to expand the latter as a finite sum.  It is this
expansion together with the Jacobian of the change of variables which
produce the variants of binomial coefficients given in (\ref{eq:22})
into our calculations.  A few facts regarding these coefficients are
necessary.  
\begin{lemma}
\label{lemma:3}
Let $j$ and $k$ be positive integers.  Then,
\begin{enumerate}
\item \label{bin1} ${\displaystyle \left(x +
      \frac{1}{x}\right)^{j-1}\left(x - \frac{1}{x}\right) =
    \sum_{m=1}^j \left[ j - 1 \atop m \right] x^{j - 2m}.  }$
\item \label{bin4} ${\displaystyle 2^k = \sum_{n=0}^k \left[{k-1 \atop
      n} \right] (k-2n). }$
\item \label{bin2} If $j$ is odd, then ${\displaystyle \frac{2^j}{j} =
    \sum_{m=0}^j \left[j-1 \atop m \right] \frac{1}{j-2m} }.$
\item \label{bin3} If $j$ is odd and $k$ is even, then 
\[
{\displaystyle
    \frac{2^{j+k}}{j(j+k)} = \sum_{m=1}^j \sum_{n=1}^k \left[j - 1
      \atop m \right] \left[k-1 \atop n
    \right] \frac{1}{j-2m} \times \frac{1}{j-2m + k -2n}.}
\]
\end{enumerate}
\end{lemma}
\begin{proof*}
  To prove \ref{bin1} we use the Binomial Theorem to expand $(x + 1/x)^{j-1}(x - 1/x)$ and collect together
  terms with like powers of $x$.  Fact \ref{bin4} follows by taking
the derivative of both sides of \ref{bin1} and setting $x=1$.
  
To prove \ref{bin2}, let $\omega$ be a path in the complex plane that does not pass
  through $z = 0$ and consider the path integral
\begin{equation}
\label{path integral}
\int_{\omega} \left(z + \frac{1}{z}\right)^{j-1} \left(1 -
\frac{1}{z^2}\right) \, dz. 
\end{equation}
If $j$ is odd then \ref{bin1} implies that the integrand consists
of even powers of $x$.  Thus the integral in \ref{path integral}
depends only on the end points of $\omega$.  Notice then that,
\[
\frac{2^j}{j} = \frac{1}{2}\int_{-2}^2 x^{j-1} \, dx =
\frac{1}{2}\int_{\omega} \left(z + \frac{1}{z}\right)^{j-1} \left(1 -
  \frac{1}{z^2}\right) \, dz,
\]
where the second equality follows from the change of variables $x
\mapsto z + 1/z$, and $\omega$ is any path in the complex plane
starting at $z=-1$ ending at $z=1$ and not passing through $z=0$.
Using \ref{bin1} and the Fundamental Theorem of Calculus we find
\[
\frac{2^j}{j} = \sum_{m=0}^j \left[ j - 1 \atop m\right]
\frac{1}{j-2m}.
\]
To prove \ref{bin3} we notice that
\begin{eqnarray*}
\frac{2^{j+k}}{j(j+k)} &=& \frac{1}{2} \int_{-2}^2 y^{k-1}
\int_{-2}^y x^{j-1} \, dx \, dy \\ &=& \frac{1}{2}\sum_{m=0}^j \left[ j-1
  \atop m \right] \int_{-2}^2 y^{k-1} \int_{-2}^{\phi_+(y)} x^{j-2m-1}
\, dx \, dy,
\end{eqnarray*}
where the second equality stems from the change of variables $x
\mapsto x + 1/x$, and $\phi_+(y) = (y + \sqrt{y^2 - 4})/2$.  Again we
use the fact that $j$ is odd to conclude that the resulting integral
is path independent.  Assuming that $k$ is even we may evaluate the
inner integral and simplify to find
\multbox
\begin{eqnarray*}
\frac{2^{j+k}}{j(j+k)} &=& \frac{1}{2}\sum_{m=0}^j \left[ j-1 \atop m \right]
\frac{1}{j - 2m} \int_{-2}^2 y^{k-1} \phi_+(y)^{j - 2m} \, dy \\
&=& \frac{1}{2}\sum_{m=0}^j \sum_{n=0}^k \left[j-1 \atop
  m\right]\left[k-1 \atop n\right] \frac{1}{j - 2m}\int_{-1}^1
y^{j+k-2m-2n-1} \, dy \\ 
&=& \frac{1}{2}\sum_{m=0}^j \sum_{n=0}^k \left[j-1 \atop
  m\right]\left[k-1 \atop n\right] \frac{1}{j - 2m} \times
\frac{2}{j-2m + k - 2n}.  
\end{eqnarray*}
\emultbox
\end{proof*}
\begin{proof}[of Lemma~\ref{lemma:2}]
We will use brute force to compute the entries of the matrix $U_{\mathbf{P}}[j,k] = \la
\gamma^{j-1}, \gamma^{k-1} \ra$ under the conditions that $j$ is odd
and $k$ is even.  

We will first evaluate $\la \gamma^{j-1}, \gamma^{k-1} \ra_{\R}$.
Define the functions $\phi_-, \phi_+: \R
\rightarrow \C$ by
\[
\phi_-(\alpha) = \frac{\alpha - \sqrt{\alpha^2 - 4}}{2} \quad \quad
\mbox{and} \quad \quad \phi_+(\alpha) = \frac{\alpha + \sqrt{\alpha^2
    - 4}}{2},
\]
where $\sqrt{\cdot}$ is a fixed branch of the square root which maps
the positive real axis to itself.  The root function of $\rho$
restricted to the real axis can be given by,
\[
\phi(\alpha) = \left\{
\begin{array}{ll}
-\phi_-(\alpha) & \quad \mbox{if} \quad \alpha < -2, \\
1 & \quad \mbox{if} -2 \leq \alpha \leq 2, \\
\phi_+(\alpha) & \quad \mbox{if} \quad \alpha > 2.
\end{array}
\right.
\]
From the definition of $\la \gamma^{j-1}, \gamma^{k-1} \ra_{\R}$ it
follows that
\begin{equation}
\label{eq:24}
\la \gamma^{j-1}, \gamma^{k-1} \ra_{\R} = 2\int_{-\infty}^{\infty}
\phi(y)^{-s} y^{k-1}
\underbrace{\int_{-\infty}^y \phi(x)^{-s} \, x^{j-1}
  \, dx}_{\mathcal{F}(y)} \, dy.  
\end{equation}The inner integral of (\ref{eq:24}) can be written as
\[
\mathcal{F}(y) = \int_{-\infty}^y \phi(x) \, x^{j-1} \, dx =
\left\{
\begin{array}{lc}
{\displaystyle \int_{-\infty}^y (-\phi_-(x))^{-s} \, x^{j-1} \, dx} &
  \quad y \leq -2, \\ & \\
{\displaystyle \mathcal{F}(-2) + \int_{-2}^y x^{j-1} \, dx } & \quad -2
  < y \leq 2, \\ & \\
{\displaystyle \mathcal{F}(2) + \int_2^y \phi_+(x)^{-s} \, x^{j-1}
  \, dx } & \quad 2 < y.
\end{array}
\right.
\]
Each of these integrals converge when $\Re(s) > j$.  When $y \leq 2$
we may use the change of variables $x \mapsto x + 1/x$ to write
\begin{eqnarray*}
\mathcal{F}(y) &=& \int_{-\infty}^{\phi_-(y)} (-x)^{-s} \left(x +
\frac{1}{x}\right)^{j-1} \left(x - \frac{1}{x}\right) \, \frac{dx}{x}
\\ &=& \sum_{m=0}^j \left[ j-1 \atop m \right]
\int_{-\infty}^{\phi_-(y)} (-x)^{-s} \, x^{j - 2m - 1} \, dx,
\end{eqnarray*}where the second equation comes from
Lemma~\ref{lemma:3}.  Similarly, when $y \geq 2$ we may write
\[
\int_2^y \phi_+(x)^{-s} x^{j-1} \, dx = \sum_{m=0}^j \left[ j-1 \atop
  m \right] \int_1^{\phi_+(y)} x^{-s} x^{j-2m-1} \, dx.
\]
Evaluating these integrals using the fact that $j$ is odd we find,
\begin{equation}
\label{eq:25}
\mathcal{F}(y) = \left\{
\begin{array}{lc}
{\displaystyle \sum_{m=0}^j \left[j-1 \atop m\right] \frac{\{-\phi_-(y)\}^{-s}
  \phi_-(y)^{j - 2m}}{-s + j - 2m}} &
  \quad y \leq -2, \\ & \\
{\displaystyle \frac{y^j + 2^j}{j} - \sum_{m=0}^j \left[j-1 \atop
  m\right] \frac{1}{-s + j - 2m}
} & \quad -2
  < y \leq 2, \\ & \\
{\displaystyle  \frac{2^{j+1}}{j} + \sum_{m=0}^j \left[j-1 \atop
  m\right] \frac{\phi_+(y)^{-s + j - 2m} - 2}{-s + j - 2m} 
} & \quad 2 < y.
\end{array}
\right.
\end{equation}Now $\la \gamma^{j-1}, \gamma^{k-1} \ra_{\R}$ is given by
\begin{equation}
\label{eq:26}
2 \Bigg\{\underbrace{\int\limits_{\infty}^{-2} \left\{-\phi_-(y)\right\}^{-s} y^{k-1} \mathcal{F}(y)
\, dy}_{\mbox{\ding{182}}} + \underbrace{\int\limits_{-2}^2 y^{k-1} \mathcal{F}(y) \, dy}_{\mbox{\ding{183}}} + \underbrace{\int\limits_2^{\infty} 
\phi_+(y)^{-s} y^{k-1} \mathcal{F}(y) \, dy}_{\mbox{\ding{184}}}
\Bigg\}.
\end{equation}
Using (\ref{eq:25}) we find 
\[
\mbox{\ding{182}} = 
\sum_{m=0}^j \left[j-1 \atop m \right] \frac{1}{-s + j -2m}\int_{-\infty}^{-2}
\left\{-\phi_-(y)\right\}^{-2s} \phi_-(y)^{j - 2m} y^{k-1} \, dy.
\]
The change of variables $y \mapsto y + 1/y$ together with
Lemma~\ref{lemma:3} yields 
\begin{equation}
\label{eq:27}
\mbox{\ding{182}} = 
\sum_{m=0}^j \sum_{n=0}^k \left[ j-1 \atop m \right] \left[ k-1
  \atop n \right] \frac{1}{-s+j-2m} \int_{-\infty}^{-1} \{-y\}^{-2s} y^{j - 2m + k - 2n -
  1} \, dy.
\end{equation}
Again substituting (\ref{eq:25}) into (\ref{eq:26}) we may write
the \ding{183} as 
\begin{equation}
\label{eq:28}
\mbox{\ding{183}} = 
\left(
\frac{2^j}{j} - \sum_{m=0}^j \left[j-1 \atop m \right]\frac{1}{-s+j-2m}
\right) \underbrace{\int_{-2}^2 y^{k-1} \, dy}_{\mbox{0 since $k$ is
  even}} + \;\; \frac{1}{j}\int_{-2}^2 y^{j+k-1}
\, dy
\end{equation}
Similarly \ding{184}  can be written as
\begin{eqnarray}
\label{eq:29}
\mbox{\ding{184}} &=& 
-2\sum_{m=0}^j \left[ {j-1 \atop m} \right]
\frac{1}{-s + j -2m} \sum_{n=0}^k \left[{k-1 \atop n} \right]
\int_1^{\infty} y^{-s+k-2n-1} \, dy \\
&& + \sum_{m=0}^j \sum_{n=0}^k \left[ {j-1 \atop m} \right] \left[{k-1
    \atop n}\right] \frac{1}{-s+j-2m}\int_1^{\infty}
y^{-2s+j-2m+k-2n-1}\, dy \nonumber \\
&& + \frac{2^{j+1}}{j}
 \sum_{n=0}^k \left[{k-1 \atop n} \right]
\int_1^{\infty} y^{-s+k-2n-1} \, dy. \nonumber
\end{eqnarray}
Evaluating the integrals in (\ref{eq:27}), (\ref{eq:28}) and
(\ref{eq:29}) we may write $\la \gamma^{j-1}, \gamma^{k-1} \ra_{\R}$
as 
\begin{eqnarray*}
\lefteqn{
\overbrace{-4\sum_{m=0}^j \sum_{n=0}^k \left[
  j-1 \atop m \right] \left[ k-1  \atop n \right] \frac{1}{s-j+2m}\times
\frac{1}{2s-j-k+2m+2n}}^{\mbox{\ding{185}}} 
+ \frac{2^{j+k+2}}{j(j+k)}} \\ && + \frac{2^{j+2}}{j}\sum_{n=0}^k \left[
  {k-1 \atop n} \right] \frac{1}{s-k+2n} + \underbrace{4\sum_{m=0}^j
\sum_{n=0}^k \left[{j-1 \atop m} \right]\left[{k-1 \atop n} \right]\frac{1}{s-j+2m}\times\frac{1}{s-k-2n}}_{\mbox{\ding{186}}}
\end{eqnarray*}Decomposing \ding{185} + \ding{186} into partial
fractions gives
\begin{eqnarray}
\lefteqn{\la \gamma^{j-1}, \gamma^{k-1} \ra_{\R} = \frac{2^{j+k+2}}{j(j+k)} +
\frac{2^{j+2}}{j} \sum_{n=0}^k \left[k-1 \atop n\right]\frac{1}{s -
  k + 2n} \nonumber} \\ && + 4 \sum_{m=0}^j \sum_{n=0}^k \left[j-1
  \atop m \right] \left[k-1 \atop n\right] \frac{1}{2m - 2n + k -
  j}\times \frac{1}{s - k +  2n} \label{eq:30}  \\  && \underbrace{- 8
\sum_{m=0}^j \sum_{n=0}^k \left[j-1 \atop m \right] \left[k-1 \atop
  n\right] \frac{1}{2m - 2n + k - j}\times \frac{1}{2 s - k -j + 2n +
  2m} }_{\mding{187}} .  \nonumber
\end{eqnarray}
We now turn our attention to $\la \gamma^{j-1}, \gamma^{k-1}
\ra_{\C}$.
\[
\la \gamma^{j-1}, \gamma^{k-1} \ra_{\C} = -2 i \int_{\C}
\phi_1(\beta)^{-2s} (\overline{\beta})^{j-1} \beta^{k-1}
\sgn \Im(\beta) \, d\lambda_(\beta).
\]
After the change of variables $\beta \mapsto \beta + 1/\beta$ we
may rewrite $\la \gamma^{j-1}, \gamma^{k-1} \ra_{\C}$ as
\begin{eqnarray*}
-i \int_{\C} \max\{|\beta|, \left|\beta^{-1}\right|\}^{-2s}
\left(\overline{\beta} + \frac{1}{\;\overline{\beta}\;}\right)^{j-1}
\left(\beta + \frac{1}{\beta}\right)^{k-1} \left|\beta -
  \frac{1}{\beta}\right|^2 \\
\hspace{1cm}  \times \sgn \Im\left(\beta +
    \frac{1}{\beta}\right) \frac{d
  \lambda_2(\beta)}{|\beta|^2}.
\end{eqnarray*}
The integrand is invariant under the map $\beta \mapsto
1/\beta$, and thus we may replace the domain of integration with $\C
\setminus D$ (recall that $D$ is the closed unit disk).  In this
domain, $\sgn(\Im(\beta + 1/\beta)) = 1$ if $\beta$ in the open upper
half plane, and is equal to $-1$ if $\beta$ is in the open lower half
plane.  After an easy simplification we may use these facts to rewrite
$\la \gamma^{j-1}, \gamma^{k-1} \ra_{\C}$ as
\[
-4i \int\limits_{H \setminus D} |\beta|^{-2s}
\left(\overline{\beta} + \frac{1}{\;\overline{\beta}\;}\right)^{j-1}
\left(\beta + \frac{1}{\beta}\right)^{k-1} \left(\overline{\beta} -
  \frac{1}{\;\overline{\beta}\;}\right) \left(\beta -
  \frac{1}{\beta}\right) \frac{d \lambda_2(\beta)}{|\beta|^2},
\]
where $H$ is the open upper half plane.  Employing Lemma~\ref{lemma:3}
we may rewrite this as
\[
-4i \sum_{m=0}^j \sum_{n=0}^k \left[j-1 \atop m \right]\left[k-1
  \atop n\right] \int\limits_{H \setminus D}
|\beta|^{-2 s} (\overline{\beta})^{j - 2m - 1} \beta^{k - 2n - 1} \,
d\lambda_2{\beta}.
\]
Switching to polar coordinates this becomes
\[
-4 i \sum_{m=0}^j \sum_{n=0}^k \left[j-1 \atop m \right]\left[k-1
  \atop n\right] \int_0^{\pi} e^{(2m - 2n + k-j)i \theta} \, d\theta
\int_{1}^{\infty} r^{-2s +k+j-2n-2m-1} \, dr.
\]
Of course, these integrals are elementary and we finally can write 
$\la \gamma^{j-1}, \gamma^{k-1} \ra_{\C}$ as
\[
8 \sum_{m=0}^j
\sum_{n=0}^k \left[j-1 \atop m \right] \left[k-1 \atop n\right]
\frac{1}{2m - 2n + k - j}\times \frac{1}{2 s - k -j + 2n + 2m}.
\]
Notice that this exactly cancels \ding{187} from (\ref{eq:30}).  Thus,
\begin{eqnarray*}
\lefteqn{\la \gamma^{j-1}, \gamma^{k-1} \ra = \overbrace{\frac{2^{j+k+2}}{j(j+k)}}^{\mding{188}} +
\overbrace{\frac{2^{j+2}}{j}}^{\mding{189}} \sum_{n=0}^k \left[k-1 \atop n\right]\frac{1}{s -
  k + 2n} \nonumber} \\ && + 4 \sum_{m=0}^j \sum_{n=0}^k \left[j-1
  \atop m \right] \left[k-1 \atop n\right] \frac{1}{2m - 2n + k -
  j}\times \frac{1}{s - k +  2n}.
\end{eqnarray*}Using Lemma~\ref{lemma:3} we may replace \ding{188} and
\ding{189} so that
\begin{eqnarray*}
&& \la \gamma^{j-1}, \gamma^{k-1} \ra = 4 \sum_{m=0}^j \sum_{n=0}^k \left[ {j-1 \atop m} \right] \left[ {k-1
    \atop n} \right] \bigg(
\frac{1}{j-2m} \times \frac{1}{j-2m+k-2n} 
\\ && \hspace{2cm} + \;\;\; \frac{1}{j-2m}\times\frac{1}{s-k+2n} \;\;\;
  + \;\;\;  \frac{1}{2m-2n+k-j}\times\frac{1}{s-k+2n} \bigg)
\end{eqnarray*}Reindexing by $m \mapsto j-2m$ and $n \mapsto k-2n$
allows us to write $\la \gamma^{j-1}, \gamma^{k-1} \ra$ as
\begin{equation}
\label{eq:23}
4 \sum_{m=-j}^j \sum_{n=-k}^k
\left[{j-1 \atop \frac{j-m}{2}} \right] \left[{k-1 \atop
    \frac{k-n}{2}} \right] \left(
\frac{1}{m(m+n)} + \frac{1}{m(s-n)} +
\frac{1}{(n-m)(s-n)}
\right).
\end{equation}Notice that since $j$ is odd, $\left[ {j-1 \atop
    j/2} \right] = 0$, so we need not worry about the denominator of the summand being identically zero.
Next we use the fact that, for instance,
\[
\left[{j-1 \atop \frac{j+m}{2}} \right] = -\left[{j-1 \atop
      \frac{j-m}{2}} \right], 
\]
to index the sums in (\ref{eq:23}) over only positive integers.  Doing
this and simplifying the resulting summand we find,
\[
\la \gamma^{j-1}, \gamma^{k-1} \ra = \sum_{n=1}^k \sum_{m=1}^j
\left[ {j-1 \atop \frac{j-m}{2}} \right] \left[ {k-1 \atop
    \frac{k-n}{2}} \right] \frac{n}{m} \left(\frac{16}{n^2 - m^2}
\right)\left(\frac{s^2}{s^2 - n^2} \right).
\]
The summand is identically zero if $m > j$ or $n > k$ thus we may
replace the upper bounds of summation with $2J$.  Likewise the summand
is identically zero unless $n$ is even and $m$ is odd.  We may thus
may reindex the sums by $m \mapsto 2m - 1$ and $n \mapsto 2n$.  Making
these changes and simplifying the resulting expression we arrive
at the formulation of $A_{\mathbf{P}}[j,k]$ given in the statement of
the lemma in the case where $k < (N+1)/2$. 

When $k = (N+1)/2$ we have 
\begin{eqnarray*}
A_{\mathbf{P}}[j,{\textstyle \frac{N+1}{2}}] &=& \int_{\R} \phi(\alpha)^{-s} \alpha^{j-1} \,
d\alpha \\
&=& \int_{-\infty}^{-2} \phi_-(\alpha)^{-s} \alpha^{j-1} \, d\alpha +
\int_{-2}^2 \alpha^{j-1} \, d\alpha \int^{\infty}_{2}
\phi_+(\alpha)^{-s} \alpha^{j-1} \, d\alpha \\
&=& \underbrace{\frac{2^{j+1}}{2}}_{\mding{190}} + 2\sum_{m=0}^j \left[ {j-1 \atop m} \right] \int_1^{\infty}
\alpha^{-s + j - 2m - 1} \, d\alpha,
\end{eqnarray*}by the same change of variables used before.  Replacing \ding{190}
with the formula given in Lemma \ref{lemma:3} and reindexing the sum
by $m \mapsto j-2m$ we find 
\begin{eqnarray*}
A_{\mathbf{P}}[j,{\textstyle \frac{N+1}{2}}] &=& 2\sum_{m=-j}^j \left[{j-1 \atop \frac{j-m}{2}}
\right] \frac{1}{s-m} + \frac{1}{m}  \\
&=& 2\sum_{m=1}^j\left[{j-1 \atop \frac{j-m}{2}} \right] \frac{1}{s-m} -
\frac{1}{s+m} + \frac{2}{m}.
\end{eqnarray*}Since $j$ is assumed to be odd we may make the substitution $j \mapsto
2j-1$.  Similarly, since the summand is only non-zero when $m$ is odd
we may reindex the sum by $m \mapsto 2m-1$.  Thus, after simplifying,
\[
A_{\mathbf{P}}[j, {\textstyle \frac{N+1}{2}}] = 4 \sum_{m=1}^J
\left[{2j-2 \atop j-m} \right] \frac{1}{2m-1} \left(\frac{s^2}{s^2 - (2m-1)^2}\right),
\]
where the change in the upper index of summation is justified since
$\left[ {2j-2 \atop j-m} \right] = 0$ if $m > j$.

Now from \ref{bin4} of Lemma~\ref{lemma:3} we have
\[
2^{2J} = \sum_{n=0}^{2J} \left[{2J-1 \atop n} \right] 2(J-n) =
2\sum_{n=0}^{2J} \left[{2J-1 \atop J-n}\right] 2n,
\]
and thus, since $J = (N+1)/2$ we have
\[
2^N = \sum_{n=0}^{N+1} \left[{N \atop \frac{N+1}{2} -n} \right] 2n.
\]
It follows that
\[
A_{\mathbf{P}}[j, {\textstyle \frac{N+1}{2}}] = 2^{-N+2} \left[{N \atop
    \frac{N+1}{2}-n} \right] 2n \sum_{m=1}^J \left[{2j-2 \atop j-m}
\right] \frac{1}{2m-1} \left(\frac{s^2}{s^2 - (2m-1)^2} \right),
\]
which after reorganization yields the formula for
$A_{\mathbf{P}}[j,(N+1)/2]$ given in the statement of the lemma.
\end{proof}

\section{The Proofs of Theorems \ref{thm:4} and \ref{thm:5}}\label{sec:proofs-theor-refthm}

The proofs of Theorem~\ref{thm:4} and especially Theorem~\ref{thm:5}
are rather technical.  In order to see past the technical details it
is worthwhile to look at the general strategy for these proofs.
Looking at $F_N(\Phi; s)$ and $H_N(\Phi; s)$ we see both integrals are
of the form 
\[
\int\limits_{\mbox{\small monic} \atop \mbox{\small coefficients}}
\hspace{-.5cm} \widetilde{\Phi}(\mathbf{b})^{-2s} \, d\lambda(\mathbf{b}),
\]
where, of course, the monic coefficient vectors we are integrating over
and the measure $\lambda$ are dependent on whether we are looking at
the real or complex moment function.  In order to evaluate this
integral we need to exploit the multiplicativity of $\Phi$ by making a
change of variables which allows us to integrate over the roots of
monic polynomials as opposed to the coefficients.  That is, we use 
maps of the sort
\[
E: \mbox{roots} \rightarrow \mbox{monic coefficients},
\]
to write something of the form
\[
\int\limits_{\mbox{\small roots}} \left\{\prod_{n=1}^N
  \phi(\gamma_n)^{-s}\right\} \, \Jac E(\boldsymbol{\gamma}) \,
d\lambda'(\boldsymbol{\gamma}), 
\]
where $\lambda'$ is the appropriate measure on the space of roots.
At this point we begin to see difficulties arising in the case of
real moment functions which do not occur for complex moment
functions.  Namely, the space of roots of real polynomials of degree
$N$ is more complicated then the space of roots of complex polynomials
of degree $N$.  To be quite explicit, the space of roots of complex
polynomials of degree $N$ is essentially just the identification space
formed from the canonical action of $S_N$ on $\C^N$.  Consequently,
\[
H_N(\Phi; s) = \frac{1}{N!} \int_{\C^N} \left\{ \prod_{n=1}^N
  \phi(\gamma_n)^{-2s} \right\} \Jac E(\boldsymbol{\gamma}) \,
d\lambda_{2N}(\boldsymbol{\gamma}).
\]
The space of roots of real polynomials of degree $N$ on the
other hand is partitioned into components determined by the possible
numbers of real and complex conjugate pairs of roots.  That is,
\[
F_N(\Phi; s) = \sum_{L,M \geq 0 \atop L+2M=N} F_{L,M}(\Phi; s),
\]
where $F_{L,M}(\Phi; s)$ is given by
\[
\frac{1}{2^M L! M!} \int\limits_{\R^L \times
  \C^M} \left\{ \prod_{\ell=1}^L \phi(\alpha_l)^{-s} \prod_{m=1}^M
    \phi(\beta)^{-s} \phi(\overline{\beta})^{-s} \right\} \,
  \Jac E(\boldsymbol{\alpha}, \boldsymbol{\beta}) \,
  d\lambda_L(\boldsymbol{\alpha}) \, d\lambda_{2N}(\boldsymbol{\beta}),
  \]
That is, $F_{L,M}(\Phi; s)$ measures the contribution to $F_N(\Phi; s)$
of polynomials with $L$ real roots and $M$ complex conjugate
pairs of roots.  In the case of $H_N(\Phi; s)$ we have chosen to
integrate over all root vectors in $\C^N$ instead of $\C^N / S^N$.
The $1/N!$ term in front of the integral compensates for the fact
that almost every polynomial gets counted $N!$ times by doing this.
Similarly the $1/(2^M L! M!)$ allows us to integrate over vectors of
roots in the expression for $F_N(\Phi; s)$.  

At this point it makes sense to resolve the ambiguity with the maps
represented by $E$:  Let
$E_N$ to denote the change of variables from vectors of complex
roots to vectors of complex monic coefficient vectors, and let
$E_{L,M}$ to denote the change of variables from vectors of $L$ real
roots and $M$ pairs of complex conjugate pairs of roots to real monic
coefficient vectors.  The Jacobians of $E_N$ and $E_{L,M}$ are 
related to the Vandermonde determinant.  To be explicit, given
$\boldsymbol{\gamma} \in \C^N$ let $V^{\boldsymbol{\gamma}}$ to be the
$N \times N$ matrix whose $j,k$ entry is given by
$V^{\boldsymbol{\gamma}}[j,k] = \gamma_j^{k-1}$ (that is
$V^{\boldsymbol{\gamma}}$ is the $N \times N$ Vandermonde matrix in
the variables $\gamma_1, \gamma_2, \ldots, \gamma_N$).  It is well
known that the Jacobian of $E_N$ at $\boldsymbol{\gamma}$ is given by
$| \det V^{\boldsymbol{\gamma}} |^2$.  Perhaps less well known, and
the content of Lemma~\ref{lemma:4}, is that the Jacobian of $E_{L,M}$ at
$(\boldsymbol{\alpha}, \boldsymbol{\beta})$ is given by $2^M | \det
V^{\boldsymbol{\alpha}, \boldsymbol{\beta}}|$ where we interpret
$(\boldsymbol{\alpha}, \boldsymbol{\beta})$ as the vector $(\overline{\beta_1}, \beta_1, \ldots,
\overline{\beta_M}, \beta_M, \alpha_1,
\ldots, \alpha_L,)$.  Using this we may write $H_N(\Phi;s)$
and $F_{L,M}(\Phi; s)$ as
\begin{equation}
\label{eq:9}
H_N(\Phi; s) = \frac{1}{N!} \int_{\C^N} \left\{\prod_{n=1}^N
  \phi(\gamma_n)^{-2s}  \right\} \left| \det V^{\boldsymbol{\gamma}}
\right|^2 \, d\lambda_{2N}(\boldsymbol{\gamma}),
\end{equation}
and
\begin{equation}
\label{eq:11}
F_{L,M}(\Phi; s) = \frac{1}{L!M!}\int\limits_{\R^L \times \C^M}
\left\{ \prod_{\ell=1}^L \phi(\alpha_l)^{-s} \prod_{m=1}^M
  \phi(\beta)^{-s} \phi(\overline{\beta})^{-s} \right\} \, \left|
\det  V^{\boldsymbol{\alpha}, \boldsymbol{\beta}} \right| \,
d\lambda_L(\boldsymbol{\alpha}) \, d\lambda_{2M}(\boldsymbol{\beta}).
\end{equation}
Here we can see another source of complexity in the evaluation of
$F_N(\Phi; s)$ which does not arise in the evaluation of $H_N(\Phi;
s)$.  Namely, in the expression for $F_{L,M}(\Phi; s)$ we have the
absolute value of a Vandermonde determinant, while in the expression
for $H_N(\Phi; s)$ we have the modulus squared of a Vandermonde
determinant.  In the latter case, we may treat $|\det
V^{\bs{\gamma}}|^2$ uniformly for each $\bs{\gamma} \in \C^N$ by
writing 
\begin{equation}
\label{eq:64}
| \det
V^{\boldsymbol{\gamma}} |^2 = \det V^{\boldsymbol{\gamma}} \cdot \overline{\det
V^{\boldsymbol{\gamma}}}.
\end{equation}
We do not have this luxury when working with $F_{L,M}(\Phi; s)$.  For
every $(\bs{\alpha}, \bs{\beta}) \in \R^L \times \C^M$, $\det
V^{\boldsymbol{\alpha}, \boldsymbol{\beta}}$ is a complex number
which is either real or purely imaginary, and we must treat each
$(\bs{\alpha}, \bs{\beta})$ differently based on whether
$\det V^{\boldsymbol{\alpha}, \boldsymbol{\beta}}$ is real and
positive, real and negative, on the positive imaginary axis or on the
negative imaginary axis.  It is at this key point that our evaluation
of $F_N(\Phi; s)$ diverges from the evaluation of $F_N(\mu; s)$ given
by Chern and Vaaler.

\subsection{The Proof of Theorem~\ref{thm:4}}\label{sec:proof-theor-refthm:4}

The evaluation of $H_N(\Phi; s)$ will give us insight into the
evaluation of $F_{L,M}(\Phi; s)$.
So far we have only used
the multiplicativity of $\Phi$ to produce a product over the roots of
a polynomial.  To exploit the appearance of this product we will
expand both determinants on the right hand side of (\ref{eq:64}) as a
sum over $S_N$, 
\[
| \det V^{\boldsymbol{\gamma}} |^2 = \Bigg( \underbrace{\sum_{\sigma \in S_N}
  \sgn(\sigma) \prod_{n=1}^N \gamma_{n}^{\sigma(n)-1}}_{\det V^{\boldsymbol{\gamma}}} \Bigg)
\Bigg( \underbrace{\sum_{\tau \in S_N}
  \sgn(\tau) \prod_{n=1}^N
  \overline{\gamma_{n}}^{\tau(n)-1}}_{\overline{\det
    V^{\boldsymbol{\gamma}}}}  \Bigg).
\]
Substituting this expression into (\ref{eq:9}), using the
linearity of the integral and combining the products gives
\[
H_N(\Phi; s) = \frac{1}{N!} \sum_{\sigma \in S_N} \sum_{\tau \in S_N}
\sgn(\sigma) \sgn(\tau) \int_{\C^N} \left\{ \prod_{n=1}^N \phi(\gamma_n)^{-2s}
  \gamma_{n}^{\sigma(n)-1} \overline{\gamma_n}^{\tau(n)-1} \right\} d\lambda_{2N}(\boldsymbol{\gamma}).
\]
Using Fubini's Theorem we may finally see the full usefulness of the
multiplicativity of $\Phi$,
\begin{equation}
\label{eq:10}
H_N(\Phi; s) = \frac{1}{N!} \sum_{\sigma \in S_N} \sum_{\tau \in S_N}
\sgn(\sigma) \sgn(\tau) \prod_{n=1}^N \la \gamma^{\sigma(n) - 1}|
\gamma^{\tau(n) - 1} \ra.
\end{equation}Of course, we must justify the use of Fubini's Theorem, but as
remarked previously since $\gamma^{\sigma(n) - 1}$ and
$\gamma^{\tau(n) - 1}$ are polynomials of degree less that $N$, they
are in $L^2(\nu_s)$ for $\Re(s) > N$.  Now, by reindexing the product by $n \mapsto \tau^{-1}(n)$ and noting that
$\sgn(\sigma) \sgn(\tau) = \sgn(\sigma \circ \tau^{-1})$ we see the
right hand side of (\ref{eq:10}) is simply the determinant of the $N
\times N$ matrix whose $j,k$ entry is given by $\la \gamma^{j-1},
\gamma^{k-1} \ra$ (see \cite[Lemma 3.1]{sinclair} for details).  That is, we have proved Theorem~\ref{thm:4} in the special
case where $\mathbf{P} = \{ 1, \gamma, \gamma^2, \ldots,
\gamma^{N-1}\}$.  In fact, the general case is trivially different
from this case by noticing that if $\mathbf{P} = \{P_1(\gamma),
P_2(\gamma), \ldots, P_N(\gamma) \}$ is any complete set of monic
polynomials and $V^{\mathbf{P},\bs{\gamma}}$ is the $N \times N$ matrix whose
$j,k$ entry is given by $V^{\mathbf{P}, \bs{\gamma}}[j,k] =
P_k(\gamma_j)$ then $\det V^{\mathbf{P}, \bs{\gamma}} = \det
V^{\bs{\gamma}}$. 

\subsection{Remarks on the Proof of Theorem \ref{thm:5}}
 Returning to $F_{L,M}(\Phi; s)$ let us see how the Pfaffian arises in
the formulation of $F_N(\Phi; s)$.  To do this we will use the
familiar formula for the Vandermonde determinant given by
\begin{equation}
\label{eq:12}
\det V^{\boldsymbol{\gamma}} = \prod_{1 \leq m < n \leq N} (\gamma_n - \gamma_m).
\end{equation}Setting $\boldsymbol{\gamma} = (\boldsymbol{\alpha},
\boldsymbol{\beta}) =( \overline{\beta_1}, \beta_1, \ldots,
\overline{\beta_M}, \beta_M, \alpha_1, \ldots, \alpha_L)$ we can use
(\ref{eq:12}) to determine whether $\det V^{\bs{\alpha}, \bs{\beta}}$
is on the positive real axis, negative real axis, positive imaginary
axis or negative imaginary axis.  In 
Lemma~\ref{lemma:5} we will demonstrate that  
\[
|\det V^{\bs{\alpha}, \bs{\beta}}| = (-i)^M \left\{
\prod_{1 \leq j < k \leq L} \sgn(\alpha_k - \alpha_j) \prod_{m=1}^M
\sgn \Im(\beta_m)
\right\} \det V^{\bs{\alpha}, \bs{\beta}}.
\]
Next we introduce the $L \times L$ antisymmetric matrix
$T^{\boldsymbol{\alpha}}$ whose $j,k$ entry is given by
$T^{\boldsymbol{\alpha}}[j,k] = \sgn(\alpha_k - \alpha_j)$.  When $L$
is even, we have the very important identity (Lemma~\ref{lemma:6})
\[
\prod_{1 \leq j < k \leq L} \sgn(\alpha_k - \alpha_j) = \Pf
T^{\boldsymbol{\alpha}}.
\]
When $N$ is even so is $L$, and thus in this situation,
\begin{equation}
\label{eq:13}
|\det V^{\bs{\alpha}, \bs{\beta}}| = (-i)^M \left\{ \prod_{m=1}^M
  \sgn \Im(\beta_m)  \right\} \Pf T^{\boldsymbol{\alpha}}  \cdot \det
V^{\bs{\alpha}, \bs{\beta}}. 
\end{equation}
In the case when $N$ (and hence $L$) is odd we will need to modify our
approach since the Pfaffian is only defined for even by even square
antisymmetric matrices.  The point of this section is to see the
general mechanism which produces the Pfaffian structure in $F_N(\Phi;
s)$; for now we will assume that we are in the easier case where $N$ is 
even.

Substituting (\ref{eq:13}) into (\ref{eq:11}) allows us to write
\begin{eqnarray*}
\lefteqn{ F_{L,M}(\Phi; s) = 
\frac{(-i)^M}{L!M!}\int\limits_{\R^L \times \C^M}
  \left\{ \prod_{\ell=1}^L \phi(\alpha_l)^{-s} \prod_{m=1}^M
  \sgn \Im(\beta_m) \; \phi(\beta)^{-s} \phi(\overline{\beta})^{-s}
\right\}} \\ && \hspace{6cm} \times  \Pf T^{\bs{\alpha}} \cdot \det
V^{\boldsymbol{\alpha}, \boldsymbol{\beta}} \, d\lambda_L(\boldsymbol{\alpha}) \,
 d\lambda_{2M}(\boldsymbol{\beta}). 
 \end{eqnarray*}
At first glance this does not look to be much of an improvement
 over (\ref{eq:11}), however the Pfaffian admits an expansion as a sum
 over $S_N$ similar to that of the determinant
 (equation~\ref{eq:52}).  This
 expansion together with the Laplace expansion of $\det V^{\bs{\alpha},
   \bs{\beta}}$ (using minors which depend only on $\bs{\alpha}$ or
 $\bs{\beta}$) will allow us to separate this expression into a sum
 over complementary minors whose summand is the product of two
 integrals: one over $\R^L$ and the other over $\C^M$.  The integrals
 over $\R^L$ will evaluate to Pfaffians of $L \times L$ antisymmetric
 matrices, the entries of which are skew-symmetric bilinear forms of
 the form $\la \alpha^{j-1}, \alpha^{k-1} \ra_{\R}$.  Similarly the
 integrals over $\C^M$ evaluate to Pfaffians of $M \times M$
 antisymmetric matrices with entries of the form $\la \beta^{j-1},
 \beta^{k-1} \ra_{\C}$.  The evaluation of these integrals is again
 dependent on Fubini's Theorem.  The combinatorics necessary to reduce
 the resulting sum of products of Pfaffians of antisymmetric matrices
 to the formula given in Theorem \ref{thm:5} is achieved using a
 combinatorial formula for the Pfaffian of a sum of antisymmetric
 matrices.  That is, the Pfaffian of a sum can be written as the sum
 of a product of Pfaffians.

\section{The Proof of Theorem \ref{thm:5}}
The proof of Theorem~\ref{thm:5} relies on several technical lemmas.
In order to clearly see the chain of reasoning used to prove
Theorem~\ref{thm:5} we will defer the proofs of these technical lemmas
until later.  First though, we must introduce some definitions and
notation.

\subsection{Definitions and Notation}
For each $K \leq N$ we define $\mf{I}_K^N$ to be the set of increasing
functions from $\{1,2,\ldots, K\}$ to $\{1,2\ldots, N\}$.  That is,
\[
\mf{I}_K^N = \big\{ \{1,2,\ldots,K\}
\stackrel{\mathfrak{t}}{\longrightarrow} \{1,2,\ldots, N\} \;\; : \;\;
\mathfrak{t}(1) < \mathfrak{t}(2) < \cdots < \mathfrak{t}(K) \big\}.
\]
Associated to each $\mathfrak{t} \in \mf{I}_K^N$ there exists a unique
$\mathfrak{t}' \in \mf{I}_{N-K}^N$ such that the images of $\mathfrak{t}$
and $\mathfrak{t}'$ are disjoint.  Each $\mathfrak{t} \in \mf{I}_K^N$
induces a unique permutation  $\iota_{\mf{t}} \in S_N$ given by 
\[
\iota_{\mf{t}}(n) = \left\{
\begin{array}{ll}
\mf{t}(n) & \mbox{if} \quad 1 \leq n \leq K, \\ 
\mf{t}'(n-K) & \mbox{if} \quad K < n \leq N.
\end{array}
\right.
\]
We define the {\it sign} of $\mf{t}$ by setting $\sgn(\mf{t}) =
\sgn(\iota_{\mf{t}})$.  The identity map in $\mf{I}_K^N$ is denoted by
$\mf{i}$.  To each $\mf{t}$ we associate the subset of the
symmetric group given by
\[
S_N(\mf{t}) = \{ \tau \in S_N : \tau(k) \mbox{ is in the image of }
\mf{t} \mbox{ for } k=1,2,\ldots, K \}.
\]
For each $\tau \in S_N(\mf{t})$ define the permutations
$\sigma_{\tau} \in S_K$ and $\pi_{\tau} \in S_{N-K}$ by specifying that
\[
\begin{array}{ccll}
\sigma_{\tau}(k) & = & \mathfrak{t}^{-1}(\tau(k)) & \quad k=1,2,\ldots,K \\ 
\pi_{\tau}(\ell) & = & \mathfrak{t}'^{-1}(\tau(K + \ell)) &
\quad \ell=1,2,\ldots,N-K.
\end{array}
\]
We may use these definitions to give an alternative (and useful)
description of the sign of $\mf{t} \in \mf{I}^N_K$.
\begin{lemma}
\label{lemma:12}
For every $\mf{t} \in \mf{I}^N_K$ and $\tau \in S_N(\mf{t})$,
\[
\sgn(\mf{t}) = \frac{\sgn(\tau)}{\sgn(\sigma_{\tau}) \sgn(\pi_{\tau})}.
\]
\end{lemma}
\begin{proof*}
Clearly $\iota_{\mf{t}} \in S_N(\mf{t})$, and $\iota_{\mf{t}}^{-1} \circ \tau$ permutes $\{1, 2, \ldots, K\}$ and $\{K+1, K+2,
\ldots, N\}$ disjointly.  The action of this permutation on $\{1, 2,
\ldots, K\}$ is exactly that given by $\sigma_{\tau}$.  Similarly,
\[
\pi_{\tau}(\ell) = (\iota_{\mf{t}}^{-1} \circ \tau)(K+\ell) - K \qquad
\mbox{for} \qquad \ell = 1, 2, \ldots, M.
\]
It follows that the cycles in the cycle decomposition of
$\iota_{\mf{t}}^{-1} \circ \tau$ are in one-to-one correspondence with
the cycles in the cycle decomposition of $\sigma_{\tau}$ together with
the cycles in the cycle decomposition of $\pi_{\tau}$.  This yeilds,
\[
\sgn(\iota_{\mf{t}}^{-1} \circ \tau) = \sgn(\sigma_{\tau}) \sgn(\pi_{\tau}).
\]
In other words,
\[
\singlebox
\sgn(\mf{t}) = \sgn(\iota_{\mf{t}}) =
\frac{\sgn(\tau)}{\sgn(\sigma_{\tau}) \sgn(\pi_{\tau})}. 
\esinglebox
\]
\end{proof*}

Given an $N \times N$ matrix $W$ and $\mf{u}, \mf{t} \in \mf{I}_K^N$,
define $W_{\mf{u},\mf{t}}$ to be the $K \times K$ minor whose $j,k$
entry is given by $W_{\mf{u}, \mf{t}}[j,k] = W[\mf{u}(j),
\mf{t}(k)]$.  The complimentary minor is given by $W_{\mf{u}',
  \mf{t}'}$.  As an example of the utility of this
notation, the Laplace expansion of the determinant can be written as
\begin{equation}
\label{eq:41}
\det W = \sgn(\mf{u}) \sum_{\mf{t} \in \mf{I}_K^N} \sgn(\mf{t}) \det
W_{\mf{u}, \mf{t}} \cdot \det W_{\mf{u}', \mf{t}'},
\end{equation}
where $\mf{u}$ is any fixed element of $\mf{I}_K^N$.  We will also use
the abbreviated notation $W_{\mf{u}}$ for $W_{\mf{u}, \mf{u}}$; this is
useful notation for working with Pfaffians since if $W$ is an
antisymmetric matrix then minors of the form $W_{\mf{u}}$ are also
antisymmetric.  

Throughout this section $L$ and $M$ will be non-negative integers such that $L
+ 2M = N$.  We also set $\mathbf{P}$ to be a fixed complete family of
monic polynomials.  We will reserve $J$ for the integer part of
$(N+1)/2$, and we will set $K$ to the integer part of $(L+1)/2$ so that
$2K + 2M = 2J$.  

We will use $\bs{\alpha} \in \R^L$ for a vector of real
variables and $\bs{\beta} \in \C^M$ for a vector of non-real
complex variables.  As before, $V^{\bs{\alpha},
  \bs{\beta}}$ will represent the $N \times N$ Vandermonde matrix in
the variables $ \overline{\beta_1},
\beta_1, \ldots, \overline{\beta_M}, \beta_M, \alpha_1, \ldots, \alpha_L$ and $E_{L,M}: \R^L \times \C^M \rightarrow \R^N$ is
the map given by $E_{L,M}(\bs{\alpha}, \bs{\beta}) = \mathbf{b}$ where 
\[
x^N + \sum_{n=1}^N b_n x^{N-n} = \prod_{\ell = 1}^L (x -
\alpha_{\ell}) \prod_{m=1}^M (x - \beta_m)(x - \overline{\beta_m}).
\]
It is easily seen that almost every $\mathbf{b} \in \R^N$ corresponds
to $2^M M! L!$ preimages under the map $E_{L,M}$.  

Matrices will be denoted by capital roman letters, subscripts will be used to
define minors of a matrix, while superscripts will be used to reflect
any variables or parameters on which the entries of the matrix are
dependent.  Thus, for instance $W^{\bs{\alpha},
  \bs{\beta}}_{\mf{i},\mf{t}}$ is a minor of $W$ with entries that
depend on $\bs{\alpha}$ and $\bs{\beta}$.

There are complications in the proof of Theorem~\ref{thm:5} for odd
$N$ which do not arise in the even $N$ case.  In spite of this
disparity we will present the even and odd cases simultaneously.  Any
structures necessary for the odd $N$ case but unnecessary for the even
$N$ case will be subscripted by $\circ$.

\subsection{Steps in the Proof}
As suggested in Section~\ref{sec:proofs-theor-refthm} we will use the
change(s) of variables $E_{L,M}$.  Let $\mathcal{D}_{L,M}$ represent
the subset of $\R^N$ which consists of coefficient vectors of monic
polynomials of degree $N$ with $L$ real roots and $M$ pairs of
(non-real) complex conjugate roots.  That is, $\mathcal{D}_{L,M}$ is
the image in $\R^N$ of $E_{L,M}$.  Clearly $\R^N$ is the disjoint
union of $\mathcal{D}_{L,M}$ over all pairs of non-negative integers
with $L + 2M=N$.  Thus,
\begin{eqnarray}
F_N(\Phi; s) &=& \sum_{(L,M)} \int_{\mathcal{D}_{L,M}}
\widetilde{\Phi}(\mathbf{b})^{-s} \, d\lambda_N(\mathbf{b}) 
\nonumber
\\
&=& \sum_{(L,M)}  \frac{1}{2^M M! L!} \int\limits_{\R^L}
\int\limits_{\C^M} 
\left\{\prod_{\ell=1}^L
\phi(\alpha_{\ell})^{-s}\right\}\left\{ \prod_{m=1}^M \phi(\beta_m)^{-s}
\phi(\overline{\beta_m})^{-s}\right\}  \label{eq:31} \\ 
&& \hspace{5cm} \times \Jac E_{L,M}(\bs{\alpha}, \bs{\beta}) \;
d\lambda_L(\bs{\alpha}) \, d\lambda_{2M}(\bs{\beta}), \nonumber
\end{eqnarray}where the sum over $(L,M)$ is understood to be over all non-negative
integers $L$ and $M$ such that $L + 2M = N$.
\begin{lemma}
\label{lemma:4}
The Jacobian of $E_{L,M}$ is given by
\[
\Jac E_{L,M}(\bs{\alpha}, \bs{\beta}) = 2^M \left|\det V^{\bs{\alpha},
    \bs{\beta}}\right|. 
\]
\end{lemma}

\begin{lemma}
\label{lemma:5}
Let $\bs{\gamma} \in \C^N$ be given by
\[
\bs{\gamma} = ( \overline{\beta_1},
\beta_1, \ldots, \overline{\beta_M}, \beta_M,\alpha_1, \ldots, \alpha_L),
\]
and let $W^{\bs{\alpha}, \bs{\beta}}$ be the $N \times N$ matrix whose $j,k$ entry
is given by
\[
W^{\bs{\alpha}, \bs{\beta}}[j,k] = P_k(\gamma_j).
\]
Then, if $\mf{i} \in \mf{I}_{2M}^N$ is the identity map on $\{1,2,\ldots, 2M\}$,
\[
\left| \det V^{\bs{\alpha}, \bs{\beta}} \right| = \sum_{\mf{t}\in
  \mf{I}_{2M}^N} \sgn(\mf{t}) \left\{ \det W^{\bs{\beta}}_{\mf{i}, \mf{t}}
  (-i)^M \prod_{m=1}^M \sgn \Im(\beta_m) \right\} \Bigg\{ \det
W_{\mf{i}', \mf{t}'}^{\bs{\alpha}} \prod_{j < k} \sgn(\alpha_k -
\alpha_j) \Bigg\} ,
\]
where as suggested by the notation, the minors $W_{\mf{i},
  \mf{t}}^{\bs{\beta}}$ and $W^{\bs{\alpha}}_{\mf{i}',
  \mf{t}'}$ of $W^{\bs{\alpha}, \bs{\beta}}$ are dependent only on
$\bs{\beta}$ and $\bs{\alpha}$ respectively.
\end{lemma}
Using Lemma~\ref{lemma:4} and Lemma~\ref{lemma:5} we may rewrite
(\ref{eq:31}) as
\begin{eqnarray*}
\lefteqn{F_N(\Phi; s) = \sum_{(L,M)} \frac{1}{M! L!} \sum_{\mf{t} \in
    \mf{I}_{2M}^N} 
\sgn(\mf{t}) \int\limits_{\R^L} \int\limits_{\C^M} \left\{\prod_{\ell=1}^L
\phi(\alpha_{\ell})^{-s}\right\}\left\{ \prod_{m=1}^M \phi(\beta_m)^{-s}
\phi(\overline{\beta_m})^{-s}\right\}} \\ && \times
 \Bigg\{
\det W_{\mf{i}', \mf{t}'}^{\bs{\alpha}} \prod_{j < k}
\sgn(\alpha_k - \alpha_j) 
\Bigg\}
 \left\{ 
\det W^{\bs{\beta}}_{\mf{i}, \mf{t}} (-i)^M \prod_{m=1}^M \sgn \Im(\beta_m)
\right\} \, d\lambda_L(\bs{\alpha}) \, d\lambda_{2M}(\bs{\beta}),
\end{eqnarray*}
and Fubini's Theorem yields
\begin{eqnarray}
\lefteqn{
F_N(\Phi; s) = \sum_{(L,M)} \sum_{\mf{t} \in \mf{I}_{2M}^N} \sgn(\mf{t})
\frac{1}{L!} 
\int\limits_{\R^L} \left\{\prod_{\ell=1}^L \phi(\alpha_{\ell})^{-s}
\prod_{j<k} \sgn(\alpha_k - \alpha_j)\right\} \, \det
W^{\bs{\alpha}}_{\mf{i}', \mf{t}'} \; d\lambda_{L}(\bs{\alpha})} \nonumber \\
&& \hspace{2cm} \times \frac{(-i)^M}{M!} \int\limits_{\C^M}
\left\{\prod_{m=1}^M \phi(\beta_m)^{-s} \phi(\overline{\beta_m})^{-s} \sgn
\Im(\beta_m)\right\} \, \det W^{\bs{\beta}}_{\mf{i},\mf{t}} \;
d\lambda_{2M}(\bs{\beta}). \label{eq:32}
\end{eqnarray}
\begin{lemma}
\label{lemma:6}
Let $K$ be the integer part of $(L+1)/2$.  Define $T^{\bs{\alpha}}$ to
be the $2K \times 2K$ antisymmetric matrix whose $j,k$ entry is given by
\[
T^{\bs{\alpha}}[j,k] = \left\{
\begin{array}{ll}
\sgn(\alpha_k - \alpha_j) & \quad \mbox{if} \quad j,k < L+1, \\
\sgn(k-j)  & \quad \mbox{otherwise.}
\end{array}
\right.
\]
Then,
\[
\prod_{1 \leq j < k \leq L} \sgn(\alpha_k - \alpha_j) = \Pf T^{\bs{\alpha}}.
\]
\end{lemma}
\begin{proof}
See \cite{deBruijn}.
\end{proof}
It is worth remarking that when $L$ is even, the first condition
defining $T^{\bs{\alpha}}$ is always in force.  Since the Pfaffian is
only defined for even rank antisymmetric matrices, the second
condition is used when $L$ is odd to create a $2K \times 2K$ antisymmetric
matrix from an $L \times L$ matrix. 

Using Lemma~\ref{lemma:6} we may rewrite (\ref{eq:32}) as
\begin{eqnarray}
\lefteqn{
F_N(\Phi; s) = \sum_{(L,M)} \sum_{\mf{t} \in \mf{I}_{2M}^N} \sgn(\mf{t}) \frac{1}{L!}
\int\limits_{\R^L} \left\{\prod_{\ell=1}^L \phi(\alpha_{\ell})^{-s} \right\}
\, \Pf T^{\bs{\alpha}} \cdot \det
W^{\bs{\alpha}}_{\mf{i}', \mf{t}'} \; d\lambda_{L}(\bs{\alpha})} \nonumber \\
&& \hspace{2cm} \times \frac{(-i)^M}{M!} \int\limits_{\C^M}
\left\{\prod_{m=1}^M \phi(\beta_m)^{-s} \phi(\overline{\beta_m})^{-s} \sgn
\Im(\beta_m)\right\} \, \det W^{\bs{\beta}}_{\mf{i},\mf{t}} \;
d\lambda_{2M}(\bs{\beta}). \label{eq:34}
\end{eqnarray}
It is necessary for our calculations to replace the $\mf{t} \in
\mf{I}_{2M}^N$ with elements of $\mf{I}_{2M}^{2J}$.  Each $\mf{t} \in \mf{I}_{2M}^N$
induces a unique $\mf{t}_{\circ} \in \mf{I}_{2M}^{2J}$ by setting
$\mf{t} = \mf{t}_{\circ}$.  Notice that $\mf{t}'$ and
$\mf{t}_{\circ}'$ differ in the fact that if $N \neq 2J$ then
$\mf{t}'_{\circ}(2J-2M) = 2J$.  Clearly, $\sgn(\mf{t}_{\circ}) =
\sgn(\mf{t})$. 
\begin{lemma}
\label{lemma:7}
Let $R$ be the $2J \times 2J$ matrix whose $j,k$ entry is given by
\[
R[j,k] = \left\{
\begin{array}{ll}
\la P_j, P_k \ra_{\R} & \quad \mbox{if} \quad j,k < N+1 \\
{\displaystyle \sgn(k-j) \int_{\R} \phi(\alpha)^{-s} \, P_{\min\{j,k\}}(\alpha) \, d\alpha}  & \quad \mbox{otherwise},
\end{array}
\right.
\]
and suppose that $\mf{t} \in \mf{I}_{2M}^N$.  Then,
\[
\frac{1}{L!}
\int\limits_{\R^L} \left\{\prod_{\ell=1}^L \phi(\alpha_{\ell})^{-s} \right\}
\, \Pf T^{\bs{\alpha}} \, \det
W^{\bs{\alpha}}_{\mf{i}', \mf{t}'} \; d\lambda_{L}(\bs{\alpha}) = \Pf R_{\mf{t_{\circ}'}}.
\]
\end{lemma}
When $N$ is odd and $\mf{t} \in \mf{I}_{2M}^{N}$ then
$R_{\mf{t}'}$ is an odd by odd matrix.  The introduction of
$\mf{t}_{\circ}$ is useful since the Pfaffian of $R_{\mf{t}_{\circ}'}$
is defined.
\begin{lemma}
\label{lemma:8}
Let $C$ be the $2J \times 2J$ matrix whose $j,k$ entry is given by
\[
C[j,k] = \left\{
\begin{array}{ll}
\la P_j, P_k \ra_{\C} & \quad \mbox{if} \quad j,k < N+1 \\
0 & \quad \mbox{otherwise},
\end{array}
\right.
\]
and suppose that $\mf{t} \in \mf{I}_{2M}^{N}$.  Then,
\[
\frac{(-i)^M}{M!} \int\limits_{\C^M}
\left\{\prod_{m=1}^M \phi(\beta_m)^{-s} \phi(\overline{\beta_m})^{-s} \sgn
\Im(\beta_m)\right\} \, \det W^{\bs{\beta}}_{\mf{i},\mf{t}} \;
d\lambda_{2M}(\bs{\beta}) = \Pf C_{\mf{t_{\circ}}}.
\]
\end{lemma}
Using Lemma~\ref{lemma:7} and Lemma~\ref{lemma:8} we may rewrite
(\ref{eq:34}) as
\begin{equation}
\label{eq:35}
F_N(\Phi; s) = 
\sum_{(L,M)} \sum_{\mf{t} \in \mf{I}_{2M}^N} \sgn(\mf{t}_{\circ}) \Pf R_{\mf{t}_{\circ}'} \cdot \Pf C_{\mf{t}_{\circ}}.
\end{equation}
If $\mf{u} \in \mf{I}_{2M}^{2J}$ then either $2J$ is in the image of
$\mf{u}$ or $2J$ is in the image of $\mf{u}'$.  Notice that if $2J$ is
in the image of $\mf{u}$ then $\Pf C_{\mf{u}} = 0$.  If $2J$ is in the
image of $\mf{u}'$ then $\mf{u}'(2J-2M) = 2J$ and hence $\mf{u} =
\mf{t}_{\circ}$ for some $\mf{t} \in \mf{I}^N_{2M}$.  Thus we may
replace the sum over $\mf{I}_{2M}^N$ in (\ref{eq:35}) with a sum over
$\mf{I}_{2M}^{2J}$.  Consequently,
\begin{eqnarray}
F_N(\Phi; s) &=& 
\sum_{(L,M)} \sum_{\mf{u} \in \mf{I}_{2M}^{2J}} \sgn(\mf{u}) \Pf
R_{\mf{u}'} \cdot \Pf C_{\mf{u}} \nonumber \\
&=& \sum_{M=0}^J \sum_{\mf{u} \in \mf{I}_{2M}^{2J}} \sgn(\mf{u}) \Pf
R_{\mf{u}'} \cdot \Pf C_{\mf{u}},
\end{eqnarray}
where the second equation follows since the summand has been made to
be independent of $L$.  

The final step in the proof of Theorem \ref{thm:5} will be
establishing the following lemma.
\begin{lemma}
\label{lemma:9}
Suppose that $R$ and $C$ are antisymmetric $2J \times 2J$ matrices, and let
$U = R + C$.  Then,
\[
\Pf U = \sum_{M=0}^J \sum_{\mf{u} \in \mf{I}_{2M}^{2J}} \sgn(\mf{u})  \Pf
R_{\mf{u}'} \cdot \Pf C_{\mf{u}}.
\]
\end{lemma}
It follows that $F_N(\Phi; s) = \Pf(R + C)$.  From the
definition of $U_{\mathbf{P}}$ we see that $U_{\mathbf{P}} = R + C$,
and hence
\[
F_N(\Phi; s) = \Pf U_{\mathbf{P}}.
\]
\subsection{The Proof of Lemma~\ref{lemma:4}}
\begin{sloppypar}
Instead of computing the Jacobian of $E_{L,M}$ we will compute the
Jacobian of the map $E'_{L,M} : \R^L \times \C^M \rightarrow \R^N$
given by $E'_{L,M}(\bs{\alpha}, \bs{\beta}) = \mathbf{b}$ where
\[
x^N + \sum_{n=1}^N b_n x^{N-n} = \prod_{\ell =1}^L (x + \alpha_{\ell})
\prod_{m=1}^M (x + \beta_m) (x + \overline{\beta_m}).
\]
That is, $E'_{L,M}(\bs{\alpha}, \bs{\beta}) = E_{L,M}(-\bs{\alpha},
-\bs{\beta})$.  Clearly $\Jac E_{L,M}(\bs{\alpha},
\bs{\beta}) = \Jac E'_{L,M}(\bs{\alpha}, \bs{\beta})$.  
\end{sloppypar}

The $n$th coordinate function of $E'_{L,M}(\boldsymbol{\alpha},
\boldsymbol{\beta})$ is given by,
\[
e_n = e_{n}(\boldsymbol{\alpha}, \boldsymbol{\beta}) = e_{n}(\alpha_1,
\ldots, \alpha_L, \overline{\beta_1}, \beta_1, \ldots, \overline{\beta_M},
\beta_M),
\]
where $e_n$ is the $n$th elementary symmetric function.  We will use
the standard convention that $e_0 = 1$ and $e_n = 0$ if $n < 0$.  We also specify that if $1
\leq \ell < L$ then $e_{n, \ell} = e_{n,\ell}(\boldsymbol{\alpha},
\boldsymbol{\beta})$ is the $n$-the elementary symmetric function in
all of our variables except $\alpha_{\ell}$.  Similarly if $1 \leq m <
M$ then we define $e'_{n,m} = e'_{n,m}(\boldsymbol{\alpha},
\boldsymbol{\beta})$ to be the $n$-th elementary symmetric function in
all of our variables except $\beta_m$ and $\overline{\beta_m}$.

Using these definitions it is easy to see that for $1 \leq \ell \leq
L$, 
\[ 
\frac{\partial e_{n}}{\partial \alpha_{\ell}} =
e_{n-1,\ell}.
\] 

Setting $\beta_m = x_m + i y_m$ we may compute the partial derivatives of
$e_{n}$ with respect to $x_m$ and $y_m$.  We may categorize the
monomials in $e_{n}$ into four types: those which contain
$\beta_m$ but not $\overline{\beta_m}$, those which contain
$\overline{\beta_m}$ but not $\beta_m$, those which contain both
$\beta_m$ and $\overline{\beta_m}$, and those which contain neither
$\beta_m$ nor  $\overline{\beta_m}$.   That is,
\[ 
e_{n} = (\beta_m +
\overline{\beta_m}) \, e'_{n - 1,m} + \beta_m \overline{\beta_m}
\, e'_{n -2,m} + e'_{n,m}.
\]
Or what amounts to the same thing,
\[ 
e_{n} = 2 x_m \, e'_{n- 1,m} + (x_m^2 + y_m^2) \,
e'_{n - 2,m} + e'_{n,m}.
\]
It follows that
\[
\frac{\partial e_{n}}{\partial x_m} = 2 e'_{n -1,m} + 2 x_m e'_{n -
  2,m} \qquad \mbox{and} \qquad  \frac{\partial e_{n}}{\partial  y_m}
= 2 y_m e'_{n - 2,m}.
\] 
The Jacobian of $E_{L,M}$ is thus $|\det J|$ where the
$j$th row of $J$ is given by  
\begin{eqnarray*}
(2 e'_{j -1,1} + 2 x_1 \, e'_{j -  2,1} \qquad
2 y_1 \, e'_{j -  2,1} \quad \cdots \quad
2 e'_{j -1,M} + 2 x_M \, e'_{j -  2,M} \qquad
2 y_M \, e'_{j -  2,M}\\
 \quad
e_{j-1,1} \qquad e_{j-1,2} \quad \cdots \quad e_{j-1, L}
).
\end{eqnarray*}

Now let $I$ be the $L \times L$ identity matrix and let 
\[
B = \frac{1}{2}
\begin{pmatrix}
 1 & 1 \\
 i & -i
\end{pmatrix}.
\]
Define $A$ to be the $N \times N$ block diagonal matrix 
\[
A =
\begin{pmatrix}
   B &    &         &   &\\
     &  B &         &   &\\
     &    &  \ddots &   &\\
     &    &         & B &\\
     &    &         &   & I
\end{pmatrix},
\]
and set $J' = J A$.  The $j$-th row of $J'$ is given
by
\begin{eqnarray}
\lefteqn{
( 
 e'_{j -1,1} +   {\beta}_1 \, e'_{j -  2,1}  \qquad
 e'_{j -1,1} +  {\overline{\beta_1}} \, e'_{j -  2,1}  \quad \cdots}
\nonumber 
\\  && \hspace{3cm}
 e'_{j -1,M} +  {\beta}_M \, e'_{j -  2,M} \qquad
 e'_{j -1,M} +  {\overline{\beta}_M} \, e'_{j -  2,M}
\label{eq:36}
\\  && \hspace{6cm}
e_{j-1,1} \qquad e_{j-1,2} \quad \cdots \quad e_{j-1, L}) \nonumber
\end{eqnarray}
and it is easily seen that $|\det(J)| = 2^{M} |\det(J')|$.  

Now, let $1 \leq \ell \leq L$ and define $f_{\ell}(x)$ to be the
polynomial
\[
f_{\ell}(x) = \prod_{k=1 \atop k \neq \ell}^L (x + \alpha_k)
\prod_{m=1}^M (x + \beta_m)(x + \overline{\beta_m}) = \sum_{n=1}^{N}
e_{n-1, \ell} \, x^{N-n}. 
\]
Similarly, for $1 \leq m \leq M$ define $g_m$ and $\tilde{g}_m$ by 
\begin{eqnarray*}
g_m(x) &=& \prod_{\ell=1}^L (x + \alpha_{\ell}) \Bigg\{(x + {\beta_m})
\prod_{k=1 \atop k \neq m}^M (x + \beta_k)(x +
\overline{\beta_k})\Bigg\} \\
&=& (x + {\beta_m})\left(\sum_{n=1}^{N-1}  e'_{n-1,m} \,
  x^{N-1-n}\right) \\ 
&=& \sum_{n=1}^{N}  (e'_{n-1,m} + \beta_m e'_{n-2,m})
x^{N-n},
\end{eqnarray*}
and
\begin{eqnarray*}
\tilde{g}_m(x) &=& \prod_{\ell=1}^L (x + \alpha_{\ell}) \Bigg\{(x +
  \overline{\beta_m}) \prod_{k=1 \atop k \neq m}^M (x + \beta_k)(x +
\overline{\beta_k})\Bigg\} \\
&=& \sum_{n=1}^{N} (e'_{n-1,m} + \overline{\beta_m}
  e'_{n-2,m}) x^{N-n}.
\end{eqnarray*}Notice that the coefficient vectors of the 
$g_m, \tilde{g}_m$ and $f_{\ell}$ appear as the columns of $J'$.  This
is useful in light of the following orthogonality relations.
By construction, $f_{\ell}(-\beta_m) = f_{\ell}(-\overline{\beta_m}) =
0$ for all $1 \leq \ell \leq L$ and $1 \leq m \leq M$, and
\[
f_{\ell}(-\alpha_k) = \left\{
\begin{array}{ll}
{\displaystyle \prod_{j \neq \ell} (-\alpha_{\ell} + \alpha_j )
\prod_{m=1}^M (-\alpha_{\ell} + \beta_m)(-\alpha_{\ell} +
\overline{\beta_m})} & \quad \mbox{if} \quad k = \ell, \\
0 & \quad \mbox{otherwise}.
\end{array}
\right.
\]
Similarly, $g_m(-\alpha_{\ell}) = \tilde{g}_m(-\alpha_{\ell}) = 
g_m(-\beta_m) = \tilde{g}_m(-\overline{\beta_m}) = 0$ for
all $1 \leq \ell \leq L$ and $1 \leq m \leq M$, and 
\[
g_m(-\overline{\beta_k}) = \left\{
\begin{array}{ll}
{\displaystyle
\prod_{\ell=1}^L (-\overline{\beta_m} + \alpha_{\ell})
\Bigg(({-\overline{\beta_m} + \beta_m}) \prod_{j\neq m}
  (-\overline{\beta_m} + \beta_j)(-\overline{\beta_m} +
  \overline{\beta_j}) \Bigg)}  & \quad \mbox{if } k=m, \\
0 & \quad \mbox{otherwise},
\end{array}
\right.
\]
and
\[
\tilde{g}_m(-\beta_k) = \left\{
\begin{array}{ll}
{\displaystyle
\prod_{\ell=1}^L ({-\beta_m} + \alpha_{\ell})
\Bigg(({{-\beta_m} + \overline{\beta_m}}) \prod_{j\neq m}
  ({-\beta_m} + \beta_j)({-\beta_m} +
  \overline{\beta_j}) \Bigg)} & \quad \mbox{if } k=m, \\
0 & \quad \mbox{otherwise}.
\end{array}
\right.
\]

Now, let $D$ be the $N
\times N$ matrix given by
\[
D = 
\begin{pmatrix}
(-\overline{\beta_1})^{N-1} &
(-\overline{\beta_1})^{N-2} &
 & -\overline{\beta_1} & 1 \\
(-{\beta_1})^{N-1} &
(-{\beta_1})^{N-2} &
\cdots  & -{\beta_1} & 1 \\
& \vdots & \ddots & \vdots &  \\
(-\overline{\beta_M})^{N-1} &
(-\overline{\beta_M})^{N-2} &
\cdots & -\overline{\beta_M} & 1 \\
(-{\beta_M})^{N-1} &
(-{\beta_M})^{N-2} &
 & -{\beta_M} & 1 \\
(-\alpha_1)^{N-1} & (-\alpha_1)^{N-2} & 
 & -\alpha_1 & 1 \\
(-\alpha_2)^{N-1} & (-\alpha_2)^{N-2} & 
\cdots  & -\alpha_2 & 1 \\
& \vdots & \ddots & \vdots &  \\
(-\alpha_L)^{N-1} & (-\alpha_L)^{N-2} & 
\cdots & -\alpha_L & 1 \\
\end{pmatrix}.
\]
Clearly, $D$ is a permutation matrix times the $N \times N$
Vandermonde matrix in the variables 
\[
 -\overline{\beta_1}, -\beta_1,
\ldots, -\overline{\beta_M}, -\beta_M,-\alpha_1, -\alpha_2, \ldots,
-\alpha_L.
\]
And thus,
\[
|\det D| = \left| \det V^{\boldsymbol{\alpha},\boldsymbol{\beta}} \right|.
\]

Now, from the definitions of $D$ and $J'$ (equation~\ref{eq:36}) we
find 
\[
D J' = \begin{pmatrix}
 g_1(-\overline{\beta_1})
& \tilde{g}_1(-\overline{\beta_1}) &  & g_M(-\overline{\beta_1}) &
\tilde{g}_M(-\overline{\beta_1}) & f_1(-\overline{\beta_1}) &   & f_L(-\overline{\beta_1})  \\
 {g}_1(-\beta_1) & \tilde{g}_1(-\beta_1) & \cdots & g_M(-\beta_1) &
\tilde{g}_M(-\beta_1) &  f_1(-\beta_1) &  \cdots & f_L(-\beta_1) 
\\
 & \vdots & \ddots & \vdots &  & \vdots & \ddots
& \vdots \\
g_1(-\overline{\beta_M})  
&  \tilde{g}_1(-\overline{\beta_M}) & \cdots & g_M(-\overline{\beta_M}) &
\tilde{g}_M(-\overline{\beta_M}) & f_1(-\overline{\beta_M}) &  \cdots & f_L(-\overline{\beta_M})  \\
 g_1(-{\beta_M}) & \tilde{g}_1(-\beta_M) &  & g_M(-\beta_M) &
\tilde{g}_M(-\beta_M) & f_1(-\beta_M) &   & f_L(-\beta_M) 
\\
 g_1(-\alpha_1)
& \tilde{g}_1(-\alpha_1) &  & g_M(-\alpha_1) &
\tilde{g}_M(-\alpha_1) & f_1(-\alpha_1) &   & f_L(-\alpha_1) \\ 
 g_1(-\alpha_2)
& \tilde{g}_1(-\alpha_2) & \cdots & g_M(-\alpha_2) &
\tilde{g}_M(-\alpha_2) & f_1(-\alpha_2) & \cdots  & f_L(-\alpha_2)  \\
 & \vdots & \ddots & \vdots &  & \vdots & \ddots
& \vdots \\
 g_1(-\alpha_L)
& \tilde{g}_1(-\alpha_L) & \cdots & g_M(-\alpha_L) &
\tilde{g}_M(-\alpha_L)  & f_1(-\alpha_L) &  \cdots & f_L(-\alpha_L)
\end{pmatrix}.
\]
But from the orthogonality relations we see that this is in fact a
diagonal matrix, and
\[
|\det(DJ')| = \left|\prod_{\ell =1}^L f_{\ell}(-\alpha_{\ell}) \prod_{m=1}^M
g_1(-\overline{\beta_1}) \tilde{g}_1(-\beta_1) \right| = \left| \det
  V^{\boldsymbol{\alpha}, \boldsymbol{\beta}} \right|^2.  
\]
But this implies that $|\det J'| = \left|\det V^{\boldsymbol{\alpha},
\boldsymbol{\beta}}\right|$, and hence 
\begin{equation*}
\Jac(E_{L,M}(\boldsymbol{\alpha}, \boldsymbol{\beta})) = |\det J| =
2^M |\det J'| = 2^M \left|\det V^{\boldsymbol{\alpha},
    \boldsymbol{\beta}}\right|. 
\end{equation*}

\subsection{The Proof of Lemma~\ref{lemma:5}}

Applying (\ref{eq:12}) to $\bs{\gamma} = (\bs{\alpha}, \bs{\beta})$, we see
that 
\begin{eqnarray*}
\det \lefteqn{V^{\boldsymbol{\alpha}, \boldsymbol{\beta}} = \left\{\prod_{j <
  k}(\alpha_k - \alpha_j) \right\} \prod_{l=1}^L \prod_{m=1}^M \left|\beta_m
- \alpha_l \right|^2} \quad \quad \\
&\times& \left\{
\prod_{m < n} \left|\beta_n - \beta_m \right|^2 \left|\beta_n -
  \overline{\beta_m} \right|^2 
\right\} \prod_{m=1}^M  2 i \Im(\beta_m). 
\end{eqnarray*}
And hence,
\begin{equation}
\left|\det V^{\bs{\alpha}, \bs{\beta}}\right| = (-i)^M \left\{
\prod_{j < k} \sgn(\alpha_k - \alpha_j) \prod_{m=1}^M \sgn \Im(\beta_m)
\right\} \det V^{\bs{\alpha}, \bs{\beta}}.
\label{eq:40}
\end{equation}As in the end of Section~\ref{sec:proof-theor-refthm:4}
we may replace the monomials in the Vandermonde matrix with any
complete family of monic polynomials without changing its
determinant.  That is, $\det V^{\bs{\alpha}, \bs{\beta}} = \det
W^{\bs{\alpha}, \bs{\beta}}$.  Using the Laplace expansion of the
determinant (\ref{eq:41}) with $\mf{u} = \mf{i} \in
\mf{I}_{2M}^N$, we see that
\[
\det W^{\bs{\alpha}, \bs{\beta}} = \sum_{\mf{t} \in \mf{I}_{2M}^N}
\sgn(\mf{t}) \det W^{\bs{\alpha}, \bs{\beta}}_{\mf{i},\mf{t}} \cdot \det
W^{\bs{\alpha}, \bs{\beta}}_{\mf{i}', \mf{t}'}.
\]
Notice that the minors of the form $W^{\bs{\alpha},
  \bs{\beta}}_{\mf{i}, \mf{t}}$ consists of elements from the first
$2M$ columns of $W^{\bs{\alpha}, \bs{\beta}}$.  These columns are not
dependent on $\bs{\alpha}$ and thus we may write these minors as
$W^{\bs{\beta}}_{\mf{i}, \mf{t}}$.  Similarly we may write the minors
of the form $W^{\bs{\alpha}, \bs{\beta}}_{\mf{i}', \mf{t}'}$ as
$W^{\bs{\alpha}}_{\mf{i}', \mf{t}'}$.  It follows that
\begin{equation}
\label{eq:42}
\det V^{\bs{\alpha}, \bs{\beta}} = \sum_{\mf{t} \in \mf{I}_{2M}^N}
\sgn(\mf{t}) \det W^{\bs{\beta}}_{\mf{i},\mf{t}} \cdot \det
W^{\bs{\alpha}}_{\mf{i}', \mf{t}'},
\end{equation}
and the Lemma follows by substituting (\ref{eq:42}) into (\ref{eq:40})
and simplifying.

\subsection{The Proof of Lemma~\ref{lemma:7}}
We start by setting
\begin{equation}
\label{eq:46}
\mding{182} = 
\frac{1}{L!} \int_{\R^L}
\det W^{\boldsymbol{\alpha}}_{\mathfrak{i}', \mathfrak{t}'} \cdot
\Pf T^{\boldsymbol{\alpha}}
\left\{\prod_{\ell=1}^L \phi(\alpha_{\ell})^{-s}\right\} \,
d\lambda_L(\boldsymbol{\alpha}),
\end{equation}
where $\mf{t}$ is an element of $\mf{I}_{2M}^N$.
Expanding $\det W^{\boldsymbol{\alpha}}_{\mathfrak{i}',
\mathfrak{t}'}$ as a sum over $S_{L}$ allows us to write \ding{182} as
\begin{equation}
\label{eq:47}
\mding{182} = 
\frac{1}{L!} \sum_{\sigma \in S_L} \sgn(\sigma)
\underbrace{\int_{\R^L} \left\{ \prod_{\ell = 1}^L
    \phi(\alpha_{\ell})^{-s} \right\} \left\{\prod_{k=1}^L 
P_{\mathfrak{t}(k)}(\alpha_{\sigma(k)}) \right\} \;
\Pf T^{\boldsymbol{\alpha}} \, d\lambda_L(\boldsymbol{\alpha})}_{\mding{183}}.
\end{equation}
Recalling that for each $\sigma \in S_L$, $\Pf T^{\sigma \cdot
\boldsymbol{\alpha}} = \sgn(\sigma) \Pf T^{\boldsymbol{\alpha}}$,
we use the change of variables $\bs{\alpha} \mapsto \sigma^{-1} \cdot
\bs{\alpha}$ to write \ding{183} as
\begin{eqnarray*}
\mding{183} = 
\sgn(\sigma^{-1}) \int_{\R^L} \left\{\prod_{\ell = 1}^L
\phi(\alpha_{\ell})^{-s} \right\} \left\{ \prod_{k=1}^L
P_{\mathfrak{t}(k)}(\alpha_{k})\right\} \; \Pf T^{\boldsymbol{\alpha}}
\, d\lambda_L(\boldsymbol{\alpha}).
\end{eqnarray*}
Substituting this into (\ref{eq:47}) we see that the sum over $S_L$
exactly cancels $1/L!$.  That is,
\begin{equation}
\label{eq:49}
\mding{182} = 
\int_{\R^L} \left\{ \prod_{\ell=1}^L \phi(\alpha_{\ell})^{-s} P_{\mathfrak{t}(\ell)}(\alpha_{\ell})\right\} \;
\Pf T^{\boldsymbol{\alpha}} \, d\lambda_{L}(\boldsymbol{\alpha}).
\end{equation}
Using Lemma~\ref{lemma:13} and setting $K$ to the integer part of
$(L+1)/2$, we may write $\Pf T^{\boldsymbol{\alpha}}$ as,
\[
\Pf T^{\bs{\alpha}} = \frac{1}{K!} \sum_{\tau \in \Pi_{2K}}
\sgn(\tau) \left\{ \prod_{k=1}^K \sgn\left(\alpha_{\tau(2k)} -
    \alpha_{\tau(2k - 1)}\right) \right\}. 
\]
Substituting this into (\ref{eq:49}) we find
\begin{equation}
\label{eq:50}
\mding{182} = \frac{1}{K!} \sum_{\tau \in \Pi_{2K}} \sgn(\tau)
 \int_{\R^L} \underbrace{\left\{\prod_{\ell=1}^{L} \phi(\alpha_\ell)^{-s} 
  P_{\mathfrak{t}(\ell)}(\alpha_{\ell}) \right\} \left\{\prod_{k=1}^K
  \sgn(\alpha_{\tau(2k)} - \alpha_{\tau(2k-1)} ) \right\}
 }_{\mding{184}}   d\lambda_L(\boldsymbol{\alpha}),
\end{equation}
If $L$ is odd, then for each $\tau \in \Pi_{2K}$ there
is a $k_{\circ}$ such that $\alpha_{\tau(2 k_{\circ})} =
\alpha_{L+1}$.  If we set $\ell_{\circ} = {\tau(2 k_{\circ} -  1)}$
then we may write \ding{184} as
\begin{eqnarray*}
  \mding{184} &=& 
  \phi(\alpha_{\ell_{\circ}})^{-s}
  P_{\mf{t}'(\ell_{\circ})}(\alpha_{\ell_{\circ}}) 
  \Bigg\{\prod_{\ell = 1 \atop \ell \neq \ell_{\circ}}^L
  \phi(\alpha_{\ell})^{-s} P_{\mf{t}'(\ell)}(\alpha_{\ell}) \Bigg\}
\Bigg\{ \prod_{k=1
    \atop k \neq k_{\circ}}^K \sgn(\alpha_{\tau(2k)} -
  \alpha_{\tau(2k-1)})\Bigg\} \\
  &=& 
  \phi(\alpha_{\ell_{\circ}})^{-s}
  P_{\mf{t}'(\ell_{\circ})}(\alpha_{\ell_{\circ}}) 
  \Bigg\{
  \prod_{k=1
    \atop k \neq k_{\circ}}^K \phi(\alpha_{\tau(2k)})^{-s}
  \phi(\alpha_{\tau(2k-1)})^{-s} \\
  && \hspace{1cm} \times P_{(\mf{t}'\circ
    \tau)(2k)}(\alpha_{\tau(2k)}) P_{(\mf{t}'\circ
    \tau)(2k-1)}(\alpha_{\tau(2k-1)}) \sgn(\alpha_{\tau(2k)} -
  \alpha_{\tau(2k-1)})\Bigg\},
\end{eqnarray*}
where the second equation comes from reindexing the first product by
$\ell \mapsto \tau^{-1}(\ell)$ together with the fact that $2(K-1) =
L-1$.  Substituting this into (\ref{eq:50}) and applying Fubini's
Theorem we find
\begin{eqnarray*}
\mding{182} &=& \frac{1}{K!} \sum_{\tau \in \Pi_{2K}} \sgn(\tau)
\int_{\R} \phi(x)^{-s} P_{(\mf{t}'\circ\tau)(2k_{\circ}-1)}(x) \, dx \\
&& \hspace{1cm} \times \Bigg\{\prod_{k = 1 \atop k \neq k_{\circ}}^K \int_{\R^2}
\phi(x)^{-s} \phi(y)^{-s} P_{(\mf{t}'\circ\tau)(2k)}(y)
P_{(\mf{t}'\circ\tau)(2k-1)}(x) \sgn(y - x) \, dx \, dy\Bigg\} \\
&=& \frac{1}{K!} \sum_{\tau \in \Pi_{2K}} \sgn(\tau) 
\Bigg\{\prod_{k=1 \atop k\neq k_{\circ}}^K
\la P_{(\mf{t}' \circ \tau)(2k-1)} , P_{(\mf{t}' \circ \tau)(2k)} \ra_{\R}
\Bigg\} \int_{\R} \phi(x)^{-s} P_{(\mf{t}'\circ\tau)(2k_{\circ}-1)}(x) \, dx .
\end{eqnarray*}The use of Fubini's Theorem is justified since the
integrals in the latter expression converge if $\Re(s) > N$.
Recalling the definition of $\mf{t}'_{\circ}$ gives $(\mf{t}'_{\circ} \circ \tau)(2k_{\circ}) = 2J$, and hence
\begin{equation}
\label{eq:51}
\mding{182} = \frac{1}{K!} \sum_{\tau \in \Pi_{2K}} \sgn(\tau) \;
  R_{\mf{t}'_{\circ}}[\tau(2 k_{\circ} - 1),\tau(2 k_{\circ})]
    \prod_{k=1 \atop k \neq k_{\circ}}^K R_{\mf{t}'_{\circ}}[\tau(2 k - 1),\tau(2 k)]
\end{equation}
Similarly, when $L$ is even, \ding{182} is given by
\begin{eqnarray}
\lefteqn{\frac{1}{K!} \sum_{\tau \in \Pi_{2K}} \sgn(\tau)
\Bigg\{\prod_{k = 1}^K \int_{\R^2}
\phi(x)^{-s} \phi(y)^{-s} P_{(\mf{t}'\circ\tau)(2k)}(y)
P_{(\mf{t}'\circ\tau)(2k-1)}(x) \sgn(y - x) \, dx \, dy\Bigg\}}
\nonumber \\
&& \hspace{1cm} =  \frac{1}{K!} \sum_{\tau \in \Pi_{2K}} \sgn(\tau) 
 \Bigg\{\prod_{k=1}^K
\la P_{(\mf{t}' \circ \tau)(2k-1)} , P_{(\mf{t}' \circ \tau)(2k)} \ra_{\R}
\Bigg\}. \hspace{3.5cm} \label{eq:53}
\end{eqnarray}

Regardless if $L$ is even or odd, (\ref{eq:53}) and
(\ref{eq:51}) imply that, 
\[
\mding{182} = \frac{1}{K!} \sum_{\tau \in \Pi_{2K}} \sgn(\tau)
\prod_{k=1}^K R_{\mf{t}'_{\circ}}[\tau(2 k - 1),\tau(2 k)] = \Pf R_{\mf{t}'_{\circ}}.
\]

\subsection{The Proof of Lemma~\ref{lemma:8}}
To prove Lemma~\ref{lemma:8} we set 
\[
\mding{185} = \frac{(-i)^M}{M!} \int\limits_{\C^M}
\left\{\prod_{m=1}^M \phi(\beta_m)^{-s} \phi(\overline{\beta_m})^{-s} \sgn
\Im(\beta_m)\right\} \, \det W^{\bs{\beta}}_{\mf{i},\mf{t}} \;
d\lambda_{2M}(\bs{\beta}).
\]
From the definition of $W_{\mf{i}, \mf{t}}^{\bs{\beta}}$ we can write
\[
\det W_{\mf{i}, \mf{t}}^{\bs{\beta}} = \sum_{\tau \in S_{2M}}
\sgn(\tau) \left\{\prod_{m=1}^M P_{(\mf{t} \circ
  \tau)(2m-1)}(\overline{\beta_m}) P_{(\mf{t} \circ \tau)(2m)}(\beta_m)\right\}.
\]
Substituting this into \ding{185} we see
\begin{eqnarray*}
\mding{185} &=& \frac{1}{M!} \sum_{\tau \in S_{2M}} \sgn(\tau) (-i)^M
\int_{\C^M} \left\{ 
\prod_{m=1}^M \phi(\overline{\beta_m})^{-s} \phi(\beta_m)^{-s} \sgn
\Im(\beta_m) \right\} \\
&& \hspace{2cm} \times \left\{ 
\prod_{n=1}^M P_{(\mf{t} \circ
  \tau)(2n-1)}(\overline{\beta_n}) P_{(\mf{t} \circ \tau)(2n)}(\beta_n)
\right\} \, d\lambda_{2M}(\bs{\beta}).
\end{eqnarray*}When $\Re(s) > N$ this integral converges, and hence we
may use Fubini's Theorem to write
\begin{eqnarray*}
\mding{185} &=& \frac{1}{2^M M!} \sum_{\tau \in S_{2M}} \sgn(\tau) \Bigg\{
\prod_{m=1}^M (-2 i) \int_{\C} \phi(\overline{\beta})^{-s}
\phi(\beta)^{-s} \\
&& \hspace{4cm} \times P_{(\mf{t} \circ
  \tau)(2m-1)}(\overline{\beta}) P_{(\mf{t} \circ \tau)(2m)}(\beta)
\sgn \Im(\beta) d\lambda_2(\beta) \Bigg\} \\
&=& \frac{1}{2^M M!} \sum_{\tau \in S_{2M}} \sgn(\tau) \prod_{m=1}^M
\la P_{(\mf{t} \circ \tau)(2m-1)}, P_{(\mf{t} \circ \tau)(2m)} \ra_{\C},
\end{eqnarray*}
which is $\Pf C_{\mf{t}}$.  But, by definition, $\mf{t} =
\mf{t}_{\circ}$, and hence $\mding{185} = \Pf C_{\mf{t}_{\circ}}$ as
desired.

\subsection{The Proof of Lemma~\ref{lemma:9}}

Before proving Lemma~\ref{lemma:9} we present two alternative
formulations of the Pfaffian.  Indeed, many authors give one of these
two formulations as the definition of the Pfaffian.  

An alternative
proof of Lemma~\ref{lemma:9} is given in \cite{MR1069389}. 
\begin{lemma}
\label{lemma:13}
Let $U$ be a $2J \times 2J$ antisymmetric matrix.  
\begin{enumerate}
\item \label{item:5} Let $\Pi_{2J}$ denote the subset of $S_{2J}$ composed of those
  $\sigma$ with $\sigma(2j) > \sigma(2j - 1).$  Then,
\[
\Pf U = \frac{1}{J!} \sum_{\tau \in \Pi_{2J}} \sgn(\tau) \prod_{j=1}^J
U[\tau(2j-1), \tau(2j)].
\]
\item \label{item:6} Let $\mathbf{v}_1, \mathbf{v}_2, \ldots, \mathbf{v}_{2J}$ be the standard
  basis for $\R^{2J}$, and let $\omega$ be the 2-form given by
  ${\displaystyle \omega
  = \sum_{j < k} U[j,k] \;   \mathbf{v}_j \wedge \mathbf{v}_k}$.  Then, 
\[
\frac{1}{J!} \underbrace{\omega \wedge \omega \wedge \cdots \wedge
  \omega}_{J} = \Pf U \cdot \mathbf{v}_1 \wedge \mathbf{v}_2 \wedge
\cdots \wedge \mathbf{v}_{2J}.
\]
\end{enumerate}
\end{lemma}
\begin{proof*}
Let $G_{2J}$ be the subgroup of $S_{2J}$ generated by the transpositions $(2j-1
  \;\;\; 2j)$ for $j=1,2,\ldots,J$.   Then if $\pi \in G_{2J}$ and
$\tau \in S_{2J}$, the antisymmetry of $U$ implies that
\[
\sgn(\tau) \prod_{j=1}^J U[\tau(2j-1), \tau(2j)] = \sgn(\tau \circ \pi) \prod_{j=1}^J
U[(\tau\circ \pi)(2j-1), (\tau\circ\pi)(2j)].
\]
It follows that we may replace the sum over $S_{2J}$ in the definition
of the Pfaffian with a sum over left cosets of $G_{2J}$.  Each coset
contains $2^J$ elements and $\Pi_{2J}$ forms a complete set of coset
representatives which establishes \ref{item:5}.

To prove \ref{item:6} we write
\[
\frac{1}{J!} \underbrace{\omega \wedge \omega \wedge \cdots \wedge
  \omega}_{J} = \frac{1}{J!} \Bigg\{\sum_{k_1 < m_1} \sum_{k_2 < m_2} \cdots
\sum_{k_J < m_J} \Bigg(\prod_{j=1}^J U[k_j, m_j]\Bigg)
\cdot \bigwedge_{\ell=1}^J  (\mathbf{v}_{k_{\ell}}\wedge \mathbf{v}_{m_{\ell}})\Bigg\}.
\]
Notice that if any two of the indices of
summation are equal, then the summand on the right hand side of this
expression is identically zero. Thus we may replace the $J$-fold sum
with a single sum over $\Pi_{2J}$ to write
\[
\frac{1}{J!} \underbrace{\omega \wedge \omega \wedge \cdots \wedge 
  \omega}_{J} = \frac{1}{J!} \Bigg\{\sum_{\tau \in \Pi_{2J}}
\Bigg(\prod_{j=1}^J U[\tau(2j-1), \tau(2j)]\Bigg)
\cdot \bigwedge_{\ell=1}^J (\mathbf{v}_{\tau(2 \ell-1)}\wedge
\mathbf{v}_{\tau(2 \ell)})\Bigg\}. 
\](The big wedge notation is unambiguous here since this wedge product
is independent of order).  Then, \ref{item:6} follows from
\ref{item:5} by noting that
\[
\singlebox
\bigwedge_{\ell=1}^J (\mathbf{v}_{\tau(2
  \ell-1)}\wedge \mathbf{v}_{\tau(2 \ell)})  = \sgn(\tau) \cdot
\mathbf{v}_1 \wedge \mathbf{v}_2 \wedge \cdots \wedge
\mathbf{v}_{2J}. 
\esinglebox
\]
\end{proof*}

\begin{proof}[of Lemma~\ref{lemma:9}]
Following \ref{item:6} in Lemma~\ref{lemma:13}, define the 2-forms
$\varrho$ and $\chi$ by 
\[
\varrho = \sum_{k < m} R[k,m] \mathbf{v}_k \wedge \mathbf{v}_m \qquad
\mbox{and} \qquad \chi = \sum_{k < m} C[k,m] \mathbf{v}_k \wedge \mathbf{v}_m,
\]
where $\mathbf{v}_1, \mathbf{v}_2, \ldots, \mathbf{v}_{2J}$ is the
standard basis for $\R^{2N}$.  Then
\begin{eqnarray}
\frac{1}{J!} \bigwedge_{j=1}^J (\chi+\varrho) &=& \frac{1}{J!}
\sum_{j=0}^J {J \choose j}  \underbrace{\chi \wedge \chi \wedge
  \cdots \wedge \chi}_{j} \; \wedge \; \underbrace{\varrho \wedge
  \varrho \wedge  \cdots \wedge \varrho}_{J-j} \nonumber \\
&=& \sum_{j=0}^J \bigg\{\frac{1}{j!}\;
\underbrace{\chi \wedge \chi \wedge \cdots \wedge
  \chi}_{j}\bigg\} \wedge \bigg\{\frac{1}{(J-j)!} \underbrace{\varrho \wedge \varrho \wedge
  \cdots \wedge \varrho}_{J-j} \bigg\}. \label{eq:54} 
\end{eqnarray}
Using the Lemma~\ref{lemma:13} and the linearity and alternating
property of the wedge product it can be established that
\[
\frac{1}{j!}\; \underbrace{\chi \wedge \chi \wedge \cdots
  \wedge \chi}_{j} = \sum_{\mf{u} \in \mf{I}_{2j}^{2J}} \Pf
C_{\mf{u}} \cdot \mathbf{v}_{\mf{u}(1)} \wedge \mathbf{v}_{\mf{u}(1)}
\wedge \cdots \wedge \mathbf{v}_{\mf{u}(2j)}
\]
An analogous formula holds for $1/(J-j)! \cdot \varrho \wedge \cdots
\wedge \varrho$.  Substituting these expressions into (\ref{eq:54}) we
find
\begin{equation}
\label{eq:55}
\frac{1}{J!} \bigwedge_{j=1}^J (\chi + \varrho) = \sum_{\mf{u} \in
  \mf{I}_{2j}^{2J}} \sum_{\mf{t} \in \mf{I}_{2J-2j}^{2J}} \Pf C_{\mf{u}}
    \cdot \Pf R_{\mf{t}} \cdot \mathbf{v}_{\mf{u}(1)} \wedge \cdots
    \wedge \mathbf{v}_{\mf{u}(2j)} \wedge \mathbf{v}_{\mf{t}(1)}
    \cdots \wedge \mathbf{v}_{\mf{t}(2J - 2j)}.
\end{equation}
It is immediately clear that the summand is identically zero unless the ranges
of $\mf{u}$ and $\mf{t}$ are disjoint---that is, unless $\mf{t} =
\mf{u}'$.  Thus we may remove the sum over $\mf{I}_{2J-2j}^{2J}$.
\[
\frac{1}{J!} \bigwedge_{j=1}^J (\chi + \varrho) = \left\{\sum_{\mf{u} \in
  \mf{I}_{2j}^{2J}} \Pf C_{\mf{u}}
    \cdot \Pf R_{\mf{u}'} \cdot \mathbf{v}_{\iota_{\mf{u}}(1)} \wedge
      \mathbf{v}_{\iota_{\mf{u}}(2)} \wedge \cdots \wedge
        \mathbf{v}_{\iota_{\mf{u}}(2J)}
\right\},
\]
where we combined the two wedge products in (\ref{eq:55}) into a
single product using the definition of $\iota_{\mf{u}}$.  Now, since $\sgn(\mf{u}) =
\sgn(\iota_{\mf{u}})$,
\[
\frac{1}{J!} \bigwedge_{j=1}^J (\chi + \varrho) = \left\{\sum_{\mf{u} \in
  \mf{I}_{2j}^{2J}} \sgn(\mf{u}) \Pf C_{\mf{u}}
    \cdot \Pf R_{\mf{u}'}\right\} \cdot \mathbf{v}_1 \wedge
  \mathbf{v}_2 \wedge \cdots \wedge \mathbf{v}_{2J}.
\]
And the lemma now follows from Lemma~\ref{lemma:13}.
\end{proof}

\subsection{The Proof of Corollary~\ref{cor:3}}
Suppose that $P, Q \in \R[x]$ are either both even or both
odd.  Then, if $\phi(-\beta) = \phi(\beta)$ for every $\beta \in \C$
it is easy to verify that $\la P, Q \ra_{\R} = \la P, Q \ra_{\C} =
0$.  Notice if $P \in \R[x]$ is odd, then
\[
\int_{\R} \phi(x)^{-s} P(x) \, dx = 0.
\]
Corollary~\ref{cor:3} is a consequence of the following lemma.
\begin{lemma}
Suppose that $U$ is a $2J \times 2J$ antisymmetric matrix such that $U[j,k]
= 0$ if $(j-k) \equiv 0 \bmod 2$.  Then,
\[
\Pf U = \det A,
\]
where $A$ is the $J \times J$ matrix whose $j,k$ entry is given by
$A[j,k] = U[2j-1, 2k]$.
\end{lemma}
\begin{proof}
\begin{equation}
\label{eq:62}
\Pf U = \frac{1}{2^J J!} \sum_{\sigma \in S_{N}} \sgn(\sigma)
\prod_{j=1}^J U[\sigma(2j - 1), \sigma(2j)].
\end{equation}
Clearly the product in this expression is 0 except for permutations
$\sigma$ such that 
\begin{equation}
\label{eq:63}
\sigma(2j - 1) - \sigma(2j) \equiv 1 \bmod 2 \qquad \mbox{for} \qquad
j=1,2,\ldots, J.
\end{equation}
Let $G_{N}$ denote the subgroup of $S_{N}$ given by
$$
G_{N} = \{ \sigma \in S_{N} : (\sigma(n) - n) \equiv 0 \bmod 2,
\;\; n = 1,2,\ldots,N\}.
$$
Given $\sigma \in S_N$ satisfying (\ref{eq:63}), define $\sigma^{\ast}
\in G_N$ by 
$$
\sigma^{\ast}(2j) = \left\{
\begin{array}{cl}
\sigma(2j) & \mbox{if } \sigma(2j) \mbox{ is even},  \\
\sigma(2j - 1) & \mbox{if } \sigma(2j) \mbox{ is odd},
\end{array}
\right.  
$$
and
$$
\sigma^{\ast}(2j - 1) = \left\{
\begin{array}{cl}
\sigma(2j) & \mbox{if }\sigma(2j) \mbox{ is odd},  \\
\sigma(2j - 1) & \mbox{if }\sigma(2j) \mbox{ is even}.
\end{array}
\right.
$$
Notice that $\sigma$ and $\sigma^{\ast}$ differ only by a product
of transpositions of the form $(2j-1, 2j)$ where $j = 1,2,\ldots,J$.
Moreover, since $U$ is an antisymmetric matrix,
$$
\sgn(\sigma) \prod_{j=1}^J U[\sigma(2j - 1), \sigma(2j)] =
\sgn(\sigma^{\ast}) \prod_{j=1}^J U[\sigma^{\ast}(2j - 1),
\sigma^{\ast}(2j)].
$$
Clearly, each $\sigma^{\ast} \in G_N$ represents $2^J$ different
permutations $\sigma \in S_N$ each of which satisfy (\ref{eq:63}).  We
may thus replace the sum over $S_N$ in (\ref{eq:62}) with a sum over
$G_N$ to find
$$
\Pf(U) = \frac{1}{J!} \sum_{\sigma^{\ast} \in G_N}
\sgn(\sigma^{\ast}) \prod_{j=1}^J U[\sigma^{\ast}(2j - 1),
\sigma^{\ast}(2j)].
$$
Now, since elements of $G_N$ permute even integers and odd integers
disjointly we have $G_N$ is isomorphic to $S_J \times S_J$, and hence
$$
\Pf(U) = \frac{1}{J!} \sum_{\tau \in S_J} \sum_{\pi \in S_J}
\sgn(\tau) \sgn(\pi) \prod_{j=1}^J U[2 \tau(j) - 1, 2 \pi(j)].
$$
But, by \cite[Lemma 3.1]{sinclair} this is exactly $\det A$.
\end{proof}

\affiliationone{
   Christopher D. Sinclair\\
   Pacific Institute for the Mathematical Sciences\\
   Vancouver, British Columbia V6T 1Z2 \\
   Canada

   \email{sinclair@math.ubc.ca}
}
\affiliationtwo{~}
\affiliationthree{%
   Current address:\\
   Department of Mathematics\\
   University of Colorado at Boulder\\
   Boulder, Colorado 80309-0395\\
   USA

   \email{christopher.sinclair@colorado.edu}
}
\affiliationfour{~}
\end{document}